\documentclass{article}

\usepackage[utf8]{inputenc}

\usepackage{amsmath}
\usepackage{amsfonts}
\usepackage{amssymb}
\usepackage{mathrsfs}
\usepackage{dsfont}
\usepackage{amstext, amsopn}

\usepackage{mathtools}

\usepackage[lined,ruled, noend]{algorithm2e}

\usepackage{listings}

\usepackage{cite}
\usepackage{enumerate}
\usepackage{fancybox}

\RequirePackage{hyperref}
\PassOptionsToPackage{pdftex,bookmarksopen,bookmarksnumbered}{hyperref}
\usepackage{hyperref}

\usepackage[colorinlistoftodos]{todonotes}

\usepackage{fullpage}

\usepackage{amsthm} 
\newtheorem{thm}{Theorem}
\newtheorem{cor}[thm]{Corollary}
\newtheorem{propn}[thm]{Proposition}

\newtheorem{assumption}{Assumption}

\newtheorem{eg}{Example}
\newtheorem{remark}[thm]{Remark}

\newcommand{\figr}[1]{Figure \ref{#1}}

\newcommand{\Bcal}{{\mathcal B}}
\newcommand{\Ccal}{{\mathcal C}}

\newcommand{\Ecal}{{\mathcal E}}
\newcommand{\Fcal}{{\mathcal F}}
\newcommand{\Gcal}{{\mathcal G}}

\newcommand{\Pcal}{{\mathcal P}}

\newcommand{\Cbb}{\mathbb{C}}

\newcommand{\Ibb}{\mathbb{I}}

\newcommand{\Nbb}{\mathbb{N}}

\newcommand{\Pbb}{\mathbb{P}}

\newcommand{\Rbb}{\mathbb{R}}

\newcommand{\alphahat}{{\widehat{\alpha}}}

\newcommand{\Xbar}{{\overline{X}}}

\newcommand{\tbar}{{\overline{t}}}

\newcommand{\alphabar}{{\overline{\alpha}}}

\newcommand{\alphatilde}{{\widetilde{\alpha}}}

\newcommand{\pbf}{{\bf{p}}}

\newcommand{\paren}[1]{\left(#1\right)}
\newcommand{\ecklam}[1]{\left[#1\right]}

\renewcommand{\Re}{{\rm Re}}
\renewcommand{\equiv}{:=}

\newcommand{\Rn}{{\Rbb^n}}

\newcommand{\Rm}{{\Rbb^m}}

\newcommand{\Rp}{{\Rbb_+}}

\newcommand{\Rth}{{\Rbb^3}}
\newcommand{\Rtw}{{\Rbb^2}}

\newcommand{\ucmin}[2]{\underset{#2}{\mbox{minimize}}~#1}

\newcommand{\set}[2]{\left\{#1\,\left|\,#2\right.\right\}}
\newcommand{\map}[3]{#1:\,#2\rightarrow #3\,}
\newcommand{\mmap}[3]{#1:\,#2\rightrightarrows #3\,}

\newcommand{\ip}[2]{\left\langle #1,~ #2\right\rangle}

\newcommand{\norm}[1]{\left\|#1\right\|}

\newcommand{\cpr}[2]{\ensuremath{\mathbb{P}\left(#1\,\middle|\,#2\right)}}

\newcommand{\im}{{\rm {\scriptstyle i}}}

\DeclareMathOperator{\dd}{d}
\DeclareMathOperator{\inv}{inv}

\DeclareMathOperator{\Id}{Id}

\DeclareMathOperator{\dom}{{dom\,}}

\DeclareMathOperator{\Supp}{supp} 

\DeclareMathOperator{\Fix}{\mathsf{Fix}\,}

\newcommand{\sd}{\partial}

\title{Stochastic Algorithms for Large-Scale Composite Optimization:  the Case of Single-Shot X-FEL Imaging
\footnotetext[0]{This research was funded by the Deutsche Forschungsgemeinschaft (DFG, German Research Foundation) – Project-ID 432680300 – SFB 1456, Project C02. DRL would like to thank Haven Sturm Luke for making his gaming computer available for this research.}}
\author{
D.\ Russell Luke\thanks{Institut f\"ur Numerische und Angewandte Mathematik, Universit\"at G\"ottingen,%
\ Lotzestr.~16--18, 37083 G\"ottingen, Germany. E-mail: \texttt{r.luke@math.uni-goettingen.de}.}
\and
Steffen Schultze\thanks{Max Planck Institute for Multidisciplinary Sciences,
Department of Theoretical and Computational Biophysics, 37077, Göttingen.
	\texttt{steffen.schultze@mpinat.mpg.de}.}
\and Helmut Grubmüller\thanks{Max Planck Institute for Multidisciplinary Sciences,
Department of Theoretical and Computational Biophysics, 37077, Göttingen.
	\texttt{hgrubmu@mpinat.mpg.de}.}
}

\date{\today}

\begin{document}

\maketitle

\begin{abstract}
We apply a recently developed framework for analyzing the convergence of 
stochastic algorithms to the general problem of large-scale nonconvex
composite optimization more generally, and nonconvex likelihood maximization in particular.
Our theory is demonstrated on a stochastic gradient descent algorithm
for determining the electron density of a molecule from random samples of its scattering 
amplitude.  Numerical results on an idealized synthetic example provide a proof of concept.  
This opens the door to a broad range of algorithmic possibilities and provides a basis 
for evaluating and comparing different strategies.
While this case study is very specific, it shares a structure that transfers 
easily to many problems of current interest, particularly in machine learning.
\end{abstract}

{\small \noindent {\bfseries 2010 Mathematics Subject Classification:}
  Primary
  65C40, 
  90C06, 
  90C26;  
    Secondary
    46N30, 
    60J05, 
    49M27, 
    65K05.\\ 
  }


\section{Introduction}
We study randomized algorithms for large-scale composite optimization models of the form
\begin{equation}\tag{$\Pcal$}\label{e:comp opt}
\min_{x\in \Rn}\sum_{j=0}^M g_j(x).
\end{equation}
This problem format is the starting point for most applications in machine learning, but our
particular interest in this model comes from electronic strucure determination in X-ray imaging.
The majority of the functions $g_j$ are assumed to be smooth and nonconvex, though hard constraints
and other sharp features can be incorporated through reserving $g_0$ and the lower registers
for nonsmooth functions.  The indicator function of a closed set $A\subset \Rn$, for instance,
is defined by
\[
 g_0(x) = \begin{cases}0&\mbox{ if } x\in A\\
           +\infty&\mbox{else.}
          \end{cases}
\]
In the literature for this simple problem format, two main {\em computational} challenges are addressed:
(a) $n$ is large, and (b) $M$ is large.  We will focus on the latter challenge, though the analytical
framework we apply is not limited to this situation alone.  Of course issues like nonsmoothness and nonconvexity
of the functions $g_j$ are important considerations, but these are more consequential for the analysis
than the implementation. The challenge of having too many variables
is mainly handled by {\em coordinate descent} methods.  For a demonstration of the application of our
framework to stochastic coordinate descent see \cite{Luke23}.   The main computational challenge for the
application we have
in mind is that the number of functions
in the sum is enormous, which these days means $M\approx 10^9$.  The standard numerical
approach to handling this situation is to apply an iterative scheme using only a subset of the
functions $g_j$ at each iteration;  randomized algorithms select a subset at random.

In a series of papers \cite{HerLukStu19a, HerLukStu23a, HerLukStu23b} Luke and collaborators developed
a general convergence theory for randomized algorithms.  The approach views such algorithms as Markov chains
and convergence is in the sense of distribution with respect to the Prokhorov-L\'evy (weak) or 
Wasserstein (strong) metrics in the space of probability distributions, depending on the assumptions. 
Verifying the assumptions of the theory boils down to analyzing the properties of the randomly 
selected/generated operators at the heart of the randomized algorithm.  The application of this theory to 
single-shot X-FEL imaging was briefly described in \cite{HerLukStu23b}.  In this context, we present here a
deep-dive into a randomized method for computing the maximum likelihood estimator.  We demonstrate the
approach on simulated data for single-shot X-FEL imaging.
While this case study is very specific, it shares a structure that transfers 
easily to many problems of current interest, particularly in machine learning.
Beyond demonstrating the application of recent theory, one of the main contributions of this 
work is to add a benchmark problem instance for testing randomized splitting algorithms.

For the specific application of single-shot X-FEL we present a gradient descent algorithm for minimizing the
negative log-likelihood of a given set of outcomes conditioned on a
parameterization of a given stochastic model.  In section
\ref{s:xfel as rfi} we show 
how this is an instance of {\em random function iterations} studied in 
\cite{HerLukStu19a, HerLukStu23a, HerLukStu23b}, and present numerical results in section \ref{s:results}.  
Since our general approach reaches far beyond the specific application, we begin in section \ref{s:rco and sf}
with a general presentation of random function iterations for composite optimization.

\section{Randomized composite optimization, stochastic feasibility, and abstract convergence}\label{s:rco and sf}
The abstract algorithmic template is very simple.  We consdier the collection of subsets
of $\{0,1,\dots,M\}$ with cardinality $m<M$,
$\Ibb\equiv \set{I\in 2^{\{0,1,\dots,M\}}}{|I|=m}$, indexed by
$i = 1,2,\dots,M_m\equiv \paren{\begin{array}{c}M\\ m\end{array}}$ and construct operators
$\map{T_i}{G}{G}$ for $G\subset\Rn$, tailored to each of the subsets $I_i\in \Ibb$.
To avoid some challenging technicalities, we will assume that $T_i$ is at least continuous
on $G$ for all $i$.  For instance, in
the simplest incarnation of \eqref{e:comp opt},
all of the functions $g_i$ are smooth and $T_i$ could be the steepest
descent operator of the partial sum,
\begin{equation}\label{e:SDi} T_i \equiv \tfrac{1}{m}\sum_{j\in I_i}T_j'\quad \mbox{ where
}\quad T_j'\equiv \paren{\Id - t_j\nabla_j g_j}.
\end{equation}
This is just one of many possibilities. The main thing to observe here is that the mapping $T_i$
is built by a convex combination of the elements $T_j'$.  The regularity
of the mapping $T_i$ is determined by the calculus of {\em almost $\alpha$-firmly nonexpansive
(a$\alpha$-fne) mappings} \cite[Proposition 2.4]{LukNguTam18}  $T_j'$ whose regularity,
in turn, is determined by the regularity of the functions $g_j$.  The abstract randomized
algorithm takes the form of Algorithm \ref{algo:rfi}.
The notation $\Nbb$ in the statement of Algorithm \ref{algo:rfi} denotes the natural numbers including $0$.
The initialization of this algorithm
is stochastic in that we choose a random variable $X_0$ with distribution $\mu_0\in \mathscr{P}(\Rn)$
where $\mathscr{P}(\Rn)$ is the space of probability measures on $\Rn$;  practically, however,
most algorithms will be initialized with a single deterministic point $x_0$, that is,
the inital distribution is a delta distribution $\mu_0=\delta_{x_0}$.
\begin{algorithm}
\SetKwInOut{Output}{Initialization}
  \Output{Select a random variable $X_0$ with distribution $\mu_0\in \mathscr{P}(\Rn)$,
  and $(\xi_{k})_{k\in\Nbb}$
   an i.i.d.\ sequence with values on $\{1,2,\dots,M_m\}$, where
   $X_0$ and  $(\xi_{k})$ are independently distributed.
   Given $\map{T_i}{\Rn}{\Rn}$ for $i=1,2,\dots,M_m$.}
    \For{$k=0,1,2,\ldots$}{
            {
            \begin{equation}\label{eq:sci}
                X_{k+1}= T_{\xi_{k}}(X_{k})
            \end{equation}
            }\\
    }
  \caption{Random Function Iteration (RFI)}\label{algo:rfi}
\end{algorithm}
Convergence of the iterates $X_{k+1}$ depends not only on the properties of the
mappings $T_i$ (quasi-contractivity, for instance) but also on {\em consistency}
of the fixed points of the individual mappings $T_i$.  A {\em necessary} condition
for {\em almost sure} convergence of the sequence of random variables $(X_k)_{k\in\Nbb}$
is that the mappings $T_i$ have common fixed points, that is $\bigcap_i \Fix T_i\neq\emptyset$
\cite[Proposition 2.5]{HerLukStu23a}.
This is a very strong assumption that allows one to dispense with a stochastic analysis
entirely, deriving convergence guarantees from the available determnistic theory
\cite{HerLukStu19a}.

In the absense of common fixed points, the conventional approach is to examine the
arithmetic mean of the sequence $\Xbar_k\equiv \tfrac{1}{k}\sum_j^k X_k$.  The
sequence of arithmetic means is referred to as the {\em ergodic sequence} in the optimization
literature, the limit of this sequence, when it exists, is the {\em Ces\`aro sum}.
We take a more general approach that can
quantify not only convergence to the Ces\`aro sum, but also convergence of higher moments
of the law of the iterates.  This is
achieved by viewing Algorithm \ref{algo:rfi} as a
{\em Markov chain} whose iterates converge {\em in distribution} to an
{\em invariant measure} of the corresponding Markov operator $\Pcal$.

\subsection{Notation and Definitions}\label{s:Markov}
Our notation is standard.  The measurable sets are given by the
Borel sigma algebra on a subset  $G\subset\Rn$, denoted by $\mathcal{B}(G)$.
When the law of $X$, denoted $\mathcal{L}(X)$, satisfies $\mathcal{L}(X) = \mu$,
we write $X \sim \mu\in \mathscr{P}(G)$.  Where convenient, we will also use the notation
$\Pbb(X\in\cdot)$ to denote the law of the random variable,
where $\mathbb{P}$ is the probability measure on some underlying probability space.
The \emph{support of the probability measure}
$\mu$ is the smallest closed set $A$, for which $\mu(A)=1$
and is denoted by $\Supp \mu$.

The distance of a point $x\in \Rn$ to a set
$A\subset \Rn$ in the metric $d$ is denoted by $d(x,A)\equiv \inf_{w\in A}d(x,w)$.
The {\em projector} onto a set $A$ is denoted by $P_A$ and $P_A(x)$ is the set of
all points where $d(x,A)$ is attained.  This is empty if $A$ is open, and a singleton
if $A$ is closed and convex; generically,
$P_A$ is a (possibly empty) set-valued mapping,
for which we use the notation $\mmap{P_A}{\Rn}{\Rn}$.

The following assumptions hold throughout.
\begin{assumption}\label{ass:1}
  \begin{enumerate}[(a)]
  \item\label{item:ass1:indep} $\xi$ and  $\xi_{0},\xi_{1}, \ldots, \xi_{k}$
    are i.i.d random variables for all $k\in\Nbb$ on a probability space
    with values on $\{1, 2,\dots, M_m\}$.  The variable $X_0$ is an
    random variable with values on $\Rn$, independent from $\xi_k$.
  \item\label{item:ass1:Phi} The function
  $\map{\Phi}{\Rn\times \{1, 2,\dots, M_m\}	}{\Rn}$, $(x,i)\mapsto T_{i}x$ is measurable.
  \end{enumerate}
\end{assumption}

Let $(X_{k})_{k \in \mathbb{N}}$ be a sequence of random variables with values on $G\subset \Rn$.
In \cite{HerLukStu23a} it is shown that the sequence of random variables $(X_k)$ generated
by Algorithm \ref{algo:rfi} is a Markov chain with transition kernel  $p$ given by
\begin{equation}\label{eq:trans kernel}
  (x\in G) (A\in
  \mathcal{B}(G)) \qquad p(x,A) \equiv
  \mathbb{P}(T_{\xi}x \in A)
\end{equation}
for the measurable \emph{update function} $\map{\Phi}{G\times \{1,2,\dotsm M_m\}}{G}$ given by
$\Phi(x,i)\equiv T_{i}x$.  Recall that a Markov chain  with {\em transition kernel} $p$ satisfies
  \begin{enumerate}[(i)]
  \item $\cpr{X_{k+1} \in A}{X_{0}, X_{1}, \ldots, X_{k}} =
    \cpr{X_{k+1} \in A}{X_{k}}$;
  \item $\cpr{X_{k+1} \in A}{X_{k}} = p(X_{k},A)$
  \end{enumerate}
for all $k \in
  \mathbb{N}$ and $A \in \mathcal{B}(G)$ almost surely in probability,
  $\mathbb{P}$-a.s.  The Markov operator $\mathcal{P}$ associated with this Markov chain
is defined pointwise for a measurable function $\map{f}{G}{\mathbb{R}}$ via
\begin{align*}
  (x\in G)\qquad \mathcal{P}f(x):= \int_{G} f(y) p(x,\dd{y}),
\end{align*}
when the integral exists. With the transition kernel defined above, we have
\begin{align*}
  \mathcal{P}f(x) =  \int_{\Omega}
  f(T_{\xi(\omega)}x) \mathbb{P}(\dd{\omega}).
\end{align*}

The {\em dual} Markov operator acting on a measure $\mu\in \mathscr{P}(G)$
is indicated by action on the right by $\mathcal{P}$:
\begin{align*}
  (A \in \mathcal{B}(G))\qquad (\mathcal{P}^{*}\mu) (A):=
  (\mu\mathcal{P}) (A) := \int_{G} \mathbb{P}(T_{\xi}x \in A) \mu(\dd{x}).
\end{align*}
The distribution of the $k$'th iterate of the Markov chain generated
by Algorithm \ref{algo:rfi} is therefore easily represented as  $\mathcal{L}(X_{k}) =
\mu_0 \mathcal{P}^{k}$.

  An invariant measure of the Markov operator $\mathcal{P}$ is any distribution $\pi\in \mathscr{P}$
  that satisfies $\pi \mathcal{P} = \pi$.  The
set of all invariant probability measures is denoted by
$\inv \mathcal{P}$.
Of course in general random variables do not converge, but
distributions associated with
the sequence of random variables $(X_k)$ of Algorithm \ref{algo:rfi}, if they converge to
anything, do so to invariant measures of the associated Markov operator.
The stochastic algorithm \eqref{algo:rfi} then solves the following
{\em stochastic fixed point problem} \cite{HerLukStu23a, HerLukStu23b}:
\begin{align}
  \label{eq:stoch_fix_probl}
  \mbox{Find}\qquad \pi\in\inv\mathcal{P}.
\end{align}
When the mappings $T_i$ have common fixed points (almost surely), the problem reduces to the
  {\em stochastic feasibility} problem studied in \cite{ButnariuFlam95, Butnariu95, HerLukStu19a}:
\begin{align}\label{eq:stoch_feas_probl}
\mbox{ Find }  x^{*} \in C :=
\set{x \in G}{\mathbb{P}(x = T_{\xi}x) =1}.
\end{align}

We use the {\em Wasserstein metric} for the space of measures to metrize
convergence of the laws of the iterates of the Markov chain. Let
  \begin{equation}\label{eq:p-probabiliy measures}
       \mathscr{P}_{2}(G) = \set{\mu \in \mathscr{P}(G)}{ \exists\, x
      \in G \,:\, \int \|x-y\|^2 \mu(\dd{y}) < \infty}
  \end{equation}
  where $\|\cdot\|$ is the Euclidean norm.
  The Wasserstein $2$-metric on $\mathscr{P}_{2}(G)$, with respect to
  the Euclidean norm $\|\cdot\|$ denoted $d_{W_{2}}$, is defined by
 \begin{equation}\label{eq:Wasserstein}
  d_{W_{2}}(\mu, \nu)\equiv \paren{\inf_{\gamma\in \Ccal(\mu, \nu)}\int_{G\times G}
\|x-y\|^2\gamma(dx, dy)}^{1/2}
 \end{equation}
 where $\Ccal(\mu, \nu)$ is the set of couplings of $\mu$ and $\nu$:
   \begin{align}
    \label{eq:couplingsDef}
    \Ccal(\mu,\nu) := \set{\gamma \in \mathscr{P}(G\times G)}{ \gamma(A
      \times G) = \mu(A), \, \gamma(G\times A) = \nu(A) \quad \forall A
      \in \mathcal{B}(G)}.
  \end{align}
Since we are considering the Wasserstein $2$-metric,
convergence in this metric implies that also the second moments converge in this metric.
For more background on the analysis of sequences of measures we refer
interested readers to
\cite{Billingsley, stroock2010probability, Villani2008, Szarek2006, Hairer2021}.

The Markov operator $\mathcal{P}$
is \emph{Feller} if $\mathcal{P}f \in C_{b}(G)$ whenever $f \in
C_{b}(G)$, where $C_{b}(G)$ is the set of bounded and continuous
functions from $G$ to $\mathbb{R}$.
Although this property is central to the theory of existence of invariant measures,
we will assume existence of invariant measures.
For us, the Feller property
comes into play in characterizing the stationary points of Algorithm \ref{algo:rfi}.
It is a short exercise to show that if $T_{i}$ is continuous for all $i\in \{1,2,\dots, M_m\}$
then the Markov operator $\mathcal{P}$ is Feller (see \cite[Theorem 4.22]{Hairer2006}).
An elementary fact from the theory of Markov chains
is that, if the Markov operator
$\mathcal{P}$ is Feller and $\pi$ is a cluster
point of the sequence of measures $(\mu_k) \equiv (\mu_{0}\Pcal^k)$ with respect to
convergence in distribution then $\pi$ is an invariant probability measure \cite[Theorem 1.10]{Hairer2021}.
Moreover, the set of invariant measures of a Feller Markov operator is closed
 with respect to the topology of convergence in distribution \cite[Section 5]{Hairer2006}.
For more on this theory readers are referred to the cited works of Hairer and also \cite{Billingsley}.

\subsection{Abstract Quantitative Convergence}
We can now state the main theorem concerning  Markov operators $\Pcal$ with
transition kernel $p$ given by \eqref{eq:trans kernel} for self mappings
  $\map{T_i}{G}{G}$.
For any $\mu_0\in \mathscr{P}_2(G)$, we denote
the distributions of the iterates of Algorithm \ref{algo:rfi} by
$\mu_{k} =  \mu_0 \mathcal{P}^{k} = \mathcal{L}(X_{k})$.
It will be assumed that $\inv\mathcal{P}\neq\emptyset$.

Convergence is quantified by an implicitly defined gauge function.
Recall that $\Gcal:[0,\infty) \to [0,\infty)$ is a \textit{gauge function} if
$\Gcal$ is continuous, strictly increasing
with $\Gcal(0)=0$, and $\lim_{t\to \infty}\Gcal(t)=\infty$.
The gauge $\Gcal$  is constructed
implicitly from another
nonnegative function $\map{\theta_{\tau,\epsilon}}{[0,\infty)}{[0,\infty)}$
with parameters $\tau>0$ and $\epsilon\geq 0$ satisfying
\begin{eqnarray}\label{eq:theta_tau_eps}
(i)~ \theta_{\tau,\epsilon}(0)=0; \quad (ii)~ 0<\theta_{\tau,\epsilon}(t)<t ~\forall t\in(0,\tbar]
\mbox{ for some }\tbar>0\end{eqnarray}
and
\begin{equation}\label{eq:gauge}
 \Gcal\paren{\paren{\frac{(1+\epsilon)t^2-\paren{\theta_{\tau,\epsilon}(t)}^2}{\tau}}^{1/2}}=
 t\quad\iff\quad
 \theta_{\tau,\epsilon}(t) = \paren{(1+\epsilon)t^2 - \tau\paren{\Gcal^{-1}(t)}^2}^{1/2}
\end{equation}
for $\tau>0$ fixed.  In the next theorem the parameter $\epsilon$ quantifies the regularity of the
generating mappings $T_i$; the parameter $\tau$ is directly computed from the
another constant $\alpha$ used to characterize the regularity of $T_i$.

In preparation for the results that follow, we will require at least one of the
additional assumptions on $\theta_{\tau,\epsilon}$.
\begin{assumption}\label{ass:msr convergence}
At least one of the following holds.
 \begin{enumerate}[(a)]
\item\label{t:msr convergence, necessary sublin} $\theta_{\tau,\epsilon}$ satisfies
\begin{equation}\label{eq:theta to zero}
    \theta_{\tau,\epsilon}^{(k)}(t)\to 0\mbox{ as }k\to\infty~\forall t\in(0,\tbar),
\end{equation}
 and the sequence $(\mu_k)$ is Fej\'er monotone with respect to $\inv\mathcal{P}\cap \mathscr{P}_2(G)$, i.e.
\begin{equation}\label{eq:Fejer}
d_{W_{2,\pbf}}\paren{\mu_{k+1},\, \pi}\leq d_{W_{2,\pbf}}(\mu_k, \pi) \quad \forall k\in\Nbb,
\forall \pi\in \inv\mathcal{P}\cap\mathscr{P}_2(G);
\end{equation}
\item\label{t:msr convergence, necessary lin+}
$\theta_{\tau,\epsilon}$ satisfies
\begin{equation}\label{eq:theta summable}
    \sum_{j=1}^\infty\theta_{\tau,\epsilon}^{(j)}(t)<\infty~\forall t\in(0,\tbar)
\end{equation}
where $\theta_{\tau,\epsilon}^{(j)}$ denotes the $j$-times composition of $\theta_{\tau,\epsilon}$.
\end{enumerate}
\end{assumption}

When $\Gcal$ is simply a linear gauge this becomes
\[
\Gcal(t)=\kappa t\quad\iff\quad
\theta_{\tau, \epsilon}(t)=\paren{(1+\epsilon)-\frac{\tau}{\kappa^2}}^{1/2}t\quad (\kappa\geq \sqrt{\tfrac{\tau}{(1+\epsilon)}}).
\]
The conditions in \eqref{eq:theta_tau_eps} in this
case simplify to $\theta_{\tau, \epsilon}(t)=\gamma t$ where
\begin{equation}\label{eq:theta linear}
 0< \gamma\equiv 1+\epsilon-\frac{\tau}{\kappa^2}<1\quad\iff\quad
\sqrt{\tfrac{\tau}{(1+\epsilon)}}\leq  \kappa\leq \sqrt{\tfrac{\tau}{\epsilon}}.
\end{equation}
In other words, $\theta_{\tau, \epsilon}(t)=\gamma t$ for $\gamma<1$ satisfies
Assumption \ref{ass:msr convergence}\eqref{t:msr convergence, necessary lin+}.
The weaker Assumption \ref{ass:msr convergence}\eqref{t:msr convergence, necessary sublin}
is used to characterize sublinear convergence.

\begin{propn}[convergence rates, Theorem 2.6, \cite{HerLukStu23b}]\label{t:msr convergence}
  Let $G\subset \Rn$ be compact,
  let $\map{T_i}{G}{G}$ be continuous for all $i\in \{1,2,\dots, M_m\}$, and
  let $\xi$ and $\xi_{k}$ ($k\in\mathbb{N}$) be i.i.d. random variables taking
  values on $\{1, 2,\dots, M_m\}$.
  Define the {\em Markov transport discrepancy}
  $\map{\Psi}{\mathscr{P}_2(G)}{\mathbb{R}_+}\cup\{+\infty\}$ by
 \begin{equation}\label{eq:Psi}
\Psi(\mu)\equiv \inf_{\pi\in\inv\mathcal{P}}\inf_{\gamma\in C_*(\mu,\pi)}
\left(\int_{G\times G}
\mathbb{E}\left[\norm{(x - T_\xi x) - (y-T_\xi y)}^2\right]\ \gamma(dx, dy)\right)^{1/2}.
 \end{equation}
  Assume furthermore:
  \begin{enumerate}[(a)]
  \item \label{t:msr convergence c}  There is at least one
  $\pi \in\inv\mathcal{P}\cap \mathscr{P}_2(G)$  where $\mathcal{P}$ is the Markov operator generated
  by $T_ix$ (see \eqref{eq:trans kernel});
  \item\label{t:msr convergence a} $T_i$ satisfies
  \begin{eqnarray}\label{eq:paafne i.e.}
    &&\exists \epsilon\in [0,1), \alpha\in (0,1):\quad \forall x\in G, \forall y\in \bigcup_{\pi\in\inv\mathcal{P}}\Supp(\pi)\\
    &&\quad\mathbb{E}\ecklam{\norm{T_\xi x - T_\xi y}^2}\leq
    (1+\epsilon)d^2(x,y) -
    \tfrac{1-\alpha}{\alpha}\mathbb{E}\left[\norm{(x - T_\xi x) - (y-T_\xi y)}^2\right];
    \nonumber
  \end{eqnarray}
  \item\label{t:msr convergence b}
    $\Psi(\pi)=0\iff \pi\in\inv\mathcal{P}$, and for all $\mu\in\mathscr{P}_{2}(G)$
  \begin{eqnarray}
    d_{W_2}(\mu, \inv\Pcal)    &\le&
    \Gcal\paren{d_{\mathbb{R}}(0,\Psi(\mu))}
    = \Gcal(\Psi(\mu))\label{e:Psi msr}
  \end{eqnarray}
    with gauge $\Gcal$ given implicitly by \eqref{eq:gauge} with
    $\tau=(1-\alpha)/\alpha$ and the function $\theta_{\tau,\epsilon}$ satisfying \eqref{eq:theta_tau_eps}
    where $t_0\equiv d_{W_2}\paren{\mu_0, \inv\mathcal{P}\cap\mathscr{P}_2(G)}<\tbar$ for all $\mu_0\in \mathscr{P}_2(G)$.
\end{enumerate}
Then for any $\mu_0\in \mathscr{P}_2(G)$
the distributions $(\mu_k)$ of the iterates of Algorithm \ref{algo:rfi}
satisfy
\begin{equation}\label{eq:gauge convergence}
d_{W_2}\paren{\mu_{k+1},\inv\mathcal{P}}
\leq \theta_{\tau,\epsilon}\paren{d_{W_2}\paren{\mu_k,\inv\mathcal{P}}}
\quad \forall k \in \mathbb{N}.
\end{equation}%
If in addition $\theta_{\tau,\epsilon}$ satisfies either one of the conditions in
Assumption \ref{ass:msr convergence}, then
$\mu_k\to \pi^{\mu_0}\in\inv\mathcal{P}\cap\mathscr{P}_2(G)$ as $k\to\infty$
in the $W_2$ metric with rate $O\paren{\theta_{\tau, \epsilon}(k)}$ in the case of
Assumption \ref{ass:msr convergence}\eqref{t:msr convergence, necessary sublin} and
in the case of Assumption \ref{ass:msr convergence}\eqref{t:msr convergence, necessary lin+}
with rate
$O(s_k(t_0))$
where $s_k(t_0)\equiv\lim_{N\to\infty}\sum_{j=k}^N \theta_{\tau,\epsilon}^{(j)}(t_0)$.
\end{propn}
Condition \eqref{eq:paafne i.e.} is a slight relaxation of the condition given in Theorem 2.6 of
\cite{HerLukStu23b} (they have the inequality holding for all $y\in G$, not just $y$ on the supports of the invariant measures)
but their proof also holds for this more general statement.  Both assumptions
\eqref{t:msr convergence b} and \eqref{t:msr convergence c} of Proposition \ref{t:msr convergence}
have deep roots in fixed point theory and
variational analysis.  The context is explained in more detail in the next section.

Recall that a sequence $(\mu_k)$ 
on the metric space $(\mathscr{P}_2(G),d_{W_{2}})$
 is said to \emph{converge R-linearly} to $\pi$ with rate $c\in [0,1)$
when there exists a $\beta>0$ such that
$d_{W_{2}}(\mu_k,\pi) \le \beta c^k\quad \forall k\in \Nbb$.
The sequence $(\mu_k)$ 
is said to converge  \emph{Q-linearly} to
$\pi$ with rate $c\in [0,1)$ if there is a constant $c\in [0,1)$ such that
$d_{W_{2}}(\mu_{k+1},\pi) \le c d_{W_{2}}(\mu_{k},\pi)\quad \forall k\in \Nbb$.
Q-linear convergence is
encountered with contractive fixed point mappings, and this leads to
a priori and a posteriori error estimates on the sequence \cite[Chapter 9]{OrtegaRheinboldt70}.
\begin{cor}[Corollary 2.7, \cite{HerLukStu23b}]\label{t:msr convergence - linear}
Under the same assumptions as in Proposition \ref{t:msr convergence},
if $\Psi$ satisfies \eqref{t:msr convergence b}
with gauge $\Gcal(t)=r \cdot t$ and constant $r $
satisfying $\sqrt{\frac{1-\alpha}{\alpha(1+\epsilon)}}\leq r \
<\sqrt{\frac{1-\alpha}{\alpha\epsilon}}$,
then the sequence of iterates $(\mu_k)$ converges
R-linearly to
some $\pi^{\mu_0}\in\inv\mathcal{P}\cap\mathscr{P}_2(G)$:
\begin{equation}%
d_{W_2}\paren{\mu_{k+1},\inv\mathcal{P}}
\leq c\, d_{W_2}\paren{\mu_k,\inv\mathcal{P}}
\end{equation}%
where $c\equiv \sqrt{1+\epsilon -\paren{\tfrac{1-\alpha}{r^2\alpha}}}<1$
and $r \geq r'$ satisfies $r\geq\sqrt{(1-\alpha)/\alpha(1+\epsilon)}$.
If $\inv\mathcal{P}$ consists of a single point then convergence is
Q-linear.
\end{cor}
\noindent It is important to emphasize that R-linear convergence
permits neither a priori nor a posteriori error estimates.

In the remainder of this study we determine the classes of functions $g_j$
in \eqref{e:comp opt} that
satisfy the assumptions of Proposition \ref{t:msr convergence} for
a number of possible candidates for the mappings $T_i$.  The
concrete application of negative log-likelihood minimization in the context of
X-FEL imaging falls into the identified function classes.

\section{Regularity}
Our main results concern convergence of Markov chains under
regularity assumptions that are lifted from the generating mappings $T_i$.  In
\cite{LukNguTam18} a framework was developed for a quantitative
convergence analysis of set-valued mappings $T_{i}$ that are one-sided Lipschitz
continuous in the sense of set-valued-mappings
with Lipschitz constant slightly greater than 1.  We begin with the regularity of
$T_{i}$ and follow this through to the regularity of the resulting Markov operator.

\subsection{Almost $\alpha$-firmly nonexpansive mappings}
\label{s:aafneie}
Let $G\subset \Rn$ and let
$\map{F}{G}{\Rn}$.
 The mapping $F$ is said to be {\em pointwise almost $\alpha$-firmly
nonexpansive  at $x_0\in G$ on $G$}, abbreviated {\em pointwise a$\alpha$-fne}
whenever
\begin{eqnarray}
&&
	\exists \epsilon\in[0,1)\mbox{ and }\alpha\in (0,1):\nonumber\\
	  \label{eq:paafne}&&	\quad	\|Fx - Fx_0\|^2 \le (1+\epsilon) \|x - x_0\|^2 -
      \tfrac{1-\alpha}{\alpha}\psi(x,x_0, Fx, Fx_0)\qquad \forall x \in G,
      \label{e:paafne}
\end{eqnarray}
where the {\em transport discrepancy} $\psi$ of  $F$ at $x$ and $x_0$
is defined by
\begin{eqnarray}
&&\!\!\!\!\!\!\!\!\psi(x,x_0, Fx, Fx_0)\equiv \nonumber\\
\label{eq:psi}
&&\!\!\!\!\!\!\!\! \|Fx- x\|^2+ \|Fx_0 - x_0\|^2 + \|Fx - Fx_0\|^2 +
\|x - x_0\|^2  - \|Fx - x_0\|^2   - \|x - Fx_0\|^2.
\end{eqnarray}
When the above inequality holds for all $x_0\in G$ then $F$ is said to be
{\em a$\alpha$-fne on $G$}.
The {\em violation} is the constant
$\epsilon$ for which \eqref{eq:paafne} holds.
When $\epsilon=0$ the mapping $F$ is said to be
    {\em (pointwise) $\alpha$-firmly nonexpansive}, abbreviated {\em (pointwise) $\alpha$-fne}.

The transport discrepancy $\psi$ ties the regularity of the mapping to the geometry of the space.
A short calculation shows that, in a Euclidean space, this has
the representation
\begin{equation} \label{eq:nice ineq}
\psi(x,x_0, Fx, Fx_0) =  \|(x-Fx)-(x_0-Fx_0)\|^2.
\end{equation}

The definition of pointwise a$\alpha$-fne
mappings in Euclidean spaces appeared first in \cite{LukNguTam18} where
they are called {\em pointwise almost averaged} mappings, following
the historical development of these notions dating back to Mann, Krasnoselskii, and others
\cite{mann1953mean, krasnoselski1955, edelstein1966, BruckReich77}.
The terminology of ``averaged mappings'' comes from Baillon, Bruck and Reich \cite{BaiBruRei78}.
We do not use this nomenclature because it does not accommodate nonlinear spaces where addition
is not defined.

Condition \eqref{eq:paafne i.e.} of Proposition \ref{t:msr convergence} is the pointwise a$\alpha$-fne property
{\em in expectation} at points in the support of invariant measures.
On a closed subset $G\subset\Rn$ for a general self-mapping $\map{T_i}{G}{G}$
for $i\in \{1,2,\dots, M_m\}$, the mapping
   $\map{\Phi}{G\times \{1,2,\dots, M_m\}}{G}$ be given by $\Phi(x,i) = T_{i}x$
is said to be {\em pointwise almost $\alpha$-firmly
nonexpansive  in expectation at $x_0\in  G$} on $ G$, abbreviated
{\em pointwise a$\alpha$-fne in expectation},
  whenever
  \begin{eqnarray}\label{eq:paafne i.e.'}
&&\exists \epsilon\in [0,1), \alpha\in (0,1):\quad \forall x \in  G,\\
&&\quad\mathbb{E}\ecklam{\|\Phi(x,\xi)-\Phi(x_0,\xi)\|^2}\leq
(1+\epsilon)\|x-x_0\|^2 -
\tfrac{1-\alpha}{\alpha}\mathbb{E}\left[\psi(x,x_0, \Phi(x,\xi), \Phi(x_0,\xi))\right].
\nonumber
\end{eqnarray}
When the above inequality holds for all $x_0\in  G$ then $\Phi$ is said to be
{\em almost $\alpha$-firmly nonexpansive (a$\alpha$-fne) in expectation on $G$}.  The
expected violation is a value of $\epsilon$ for which \eqref{eq:paafne i.e.'} holds.
When the violation is $0$, the qualifier ``almost'' is dropped and the abbreviation
{\em $\alpha$-fne in expectation} is used.

\begin{eg}[steepest descent mappings of smooth functions are a$\alpha$-fne]\label{eg:sd aafne}
Let $\map{g}{\Rn}{\Rbb}$ be continuously differentiable with Lipschitz, hypomonotone gradient, that is $g$ satisfies
		\begin{subequations}\label{e:g reg}
		 \begin{eqnarray}
            ~\exists L>0:\quad
            \|\nabla g(x) - \nabla g(y)\|^2&\leq& L^2\|x-y\|^2\quad \forall x,y\in\Rn,
            \label{e:g grad-Lip}
		 \end{eqnarray}
        and
        \begin{eqnarray}
            ~\exists \tau\geq 0:-\tau\norm{x-y}^2&\leq&
            \ip{\nabla g(x)-\nabla g(y)}{x-y}\quad \forall x,y\in\Rn.
            \label{e:hypomonotone}
        \end{eqnarray}
		\end{subequations}
		Then a specialization of \cite[Proposition 13]{Luke23} to just a single block mapping establishes that
		the gradient descent mapping defined by $T_{GD}\equiv \Id - t\nabla g$
is a$\alpha$-fne on $\Rn$ with violation at most
\begin{subequations}
\begin{equation}\label{e:nabla f aalph-fne violation}
\epsilon_{GD} = \left\{2t\tau + \tfrac{t^2L^2}{\alpha}\right\}<1,
\quad\mbox{ with constant }\quad  \alpha
\end{equation}
whenever the steps $t$ satisfy
\begin{equation}\label{e:nabla f aalph-fne step}
t\in\paren{0,\frac{\alpha\sqrt{\tau^2+ L^2} - \alpha\tau}{L^2}}.
\end{equation}
\end{subequations}
\noindent  If $g$ is convex then, with step size
$t < \tfrac{2\alpha}{L}$
for $\alpha\in (0,1)$,
the gradient descent mapping $T_{GD}$ is $\alpha$-fne with
constant $\alpha$ (no violation).

Notice the interdependence of the step $t$ on the constant $\alpha$:  if
one is willing to tolerate a larger $\alpha$, then a larger step $t$ is possible;  this is at the cost, however, of
increasing the violation.  Also, the closer $\alpha$ is to $1$, the more expansive the steepest descent
operator is, as seen by \eqref{eq:paafne}.  On the other hand, short steps mean slow progress.

Now, using the calculus of a$\alpha$-fne mappings \cite[Proposition 2.4]{LukNguTam18} this yields
immediately the regularity of the steepest descent mapping of the partial sum in
\eqref{e:SDi}.  Indeed,
$T_{GD_i}\equiv \tfrac{1}{m}\sum_{j\in I_i} \Id - t_j\nabla g_j$ for $I_i\in \Ibb\equiv \set{I\in 2^{\{0,1,\dots,M\}}}{|I|=m}$
is a$\alpha$-fne on $\Rn$ with violation at most
\begin{subequations}
\begin{equation}\label{e:nabla gi aalph-fne violation}
\epsilon_{GD_i} = \tfrac{1}{m}\sum_{j\in I_i}\left\{2t_j\tau_j + \tfrac{t_j^2L_j^2}{\alphabar_i}\right\}<1,
\quad\mbox{ with constant }\quad  \alphabar_i\equiv\max_{j\in I_i}{\alpha_j}
\end{equation}
whenever the steps $t_j$ satisfy
\begin{equation}\label{e:nabla gi aalph-fne step}
t_j\in\paren{0,\frac{\alphabar_i\sqrt{\tau_j^2+ L_j^2} - \alphabar_i\tau_j}{L_j^2}}.
\end{equation}
\end{subequations}
Here $L_j$ and $\tau_j$ are the respective Lipschitz and hypomonotonicity constants of the individual
functions $g_j$ defined in \eqref{e:g reg}.

In practice, it can also be expensive to draw a new sample $I_i\in \Ibb$ because this involves
moving data into and out of memory on a computer.  We therefore extend the basic partial gradient descent
approach above to allow the option of taking several steepest descent steps on each single sample $I_i$.
Here we have
\begin{equation}\label{e:SDq}
  T_{i}\equiv T_{GD_i}^q = \paren{\tfrac{1}{m}\sum_{j\in I_i} \Id - t_j\nabla g_j}^q
\end{equation}
where $q\geq 1$ is some  fixed number of
times the steepest descent operator $T_{GD_i}$ is applied.  Again using the calculus of a$\alpha$-fne mappings
\cite[Proposition 2.4]{LukNguTam18} yields the regularity estimate that $T_i$
is a$\alpha$-fne on $\Rn$ with violation at most
\begin{equation}\label{e:nabla giq aalph-fne violation}
\epsilon_{i} = \paren{1+\epsilon_{GD_i}}^q-1,
\quad\mbox{ with constant at most}\quad  \alpha_i\equiv \frac{q}{q-1+1/\alphabar_i}
\end{equation}
whenever the steps $t_j$ satisfy \eqref{e:nabla gi aalph-fne step}.
\end{eg}

\begin{eg}[incorporating nonsmooth functions]\label{eg:fb aafne}
 We mentioned that the first functions in the sum \eqref{e:comp opt} are
 reserved for constraints and other nonsmooth functions.  Suppose, that these
 are at least lower semicontinuous (lsc)
 and {\em prox-friendly}, meaning that the resolvent of the subdifferential has a closed form.
 By {\em subdifferential}, denoted $\sd g_j$, we mean the collection of all general/limiting
 {\em  subgradients} \cite[Definition 8.3]{VA}.
 Then one would use the {\em resolvents}  of the subdifferentials of these functions:
 $T'_j= J_{\sd g_j, t_j}\equiv \paren{\Id + t_j\sd g_j}^{-1}$.  For
 $g_j$ convex this is $\alpha$-fne with $\alpha=1/2$ \cite{Moreau65}.  More generally, we will assume
 that  $g_j$ is subdifferentially regular (the limiting subdifferential is nonempty on $\dom g_j$)
 with subdifferentials satisfying
    \begin{eqnarray}\exists \tau_{j}\geq 0: &&\forall x, u, v\in G\subset\Rn,
    ~\forall z\in t_j\sd g_j(u; x), w\in t_j\sd g_j(v; x), \nonumber \\
        &&-\tfrac{\tau_{j}}{2}\norm{(u+z)
        -(v+w)}^2 \leq \ip{z-w}{u-v}.
        \label{e:sd f submon}
    \end{eqnarray}%
    Then by \cite[Proposition 2.3]{LukNguTam18} the resolvent
    $J_{\partial g_j, t_j}$ is a$\alpha$-fne with constant
    $\alpha_{j}=1/2$ and violation $\tau_j$ on $G$.

    The regularity characterized by \eqref{e:sd f submon} was introduced in
    \cite{LukNguTam18} and was shown in \cite{LukTam23} to be a generalization of
    {\em hypomonotonicity} \cite{Spin81, PolRock96a,PolRockThib00};  in particular,
    {\em weakly convex} functions (functions that, upon addition by a quadratic, can be made convex)
    satisfy \eqref{e:sd f submon}, but the converse does not hold.

 There are a number of ways to combine steepest descent steps with resolvents;  we focus on
 the following generalized forward-backward mapping.
 Let $\map{g_j}{\Rn}{\Rbb\cup\{+\infty\}}$ be lsc,  prox-friendly and satisfy
 \eqref{e:sd f submon} for $j=0,1,2,\dots,s$ and satisfy \eqref{e:g reg} for $j=s+1, \dots,M$.
 Then
 \begin{equation}\label{e:fbi}
  T_i\equiv \paren{\paren{\prod_{j\in I_i, j\leq s} J_{\partial g_j, t_j}}\circ
  \paren{\tfrac{1}{m}\sum_{j\in I_i, j> s}\Id - t_j\nabla g_j}^q}^r
 \end{equation}
where $q,r\geq 1$ allow for repeated application of the resepctive operators, and the product
$\prod_{j\in I_i, j\leq s} J_{\partial g_j, t_j}=1$ when there are no indexes $j\in I_i$ with
$j\leq s$.
Application of the calculus of a$\alpha$-fne mappings establishes that $T_i$ given by
\eqref{e:fbi} is a$\alpha$-fne on $G$ with violation
\begin{subequations}\label{e:nabla giqr aalph-fne}
\begin{equation}\label{e:nabla giqr aalph-fne violation}
\epsilon_{i} = \paren{\paren{\prod_{j\in I_i, j\leq s}\paren{1+\tau_j}}\paren{1+\epsilon_{GD_i}}^q}^r-1,
\end{equation}
and with constant at most
\begin{equation}\label{e:nabla giqr aalph-fne const}
\alphahat_i\equiv \frac{r}{r-1+1/\alphatilde_i}\quad{where}\quad
\alphatilde_i\equiv \frac{|I_i\setminus \{j>s\}|+1}{|I_i\setminus \{j>s\}| + \max\{1/2, \alpha_i\}}
\end{equation}
\end{subequations}
for $\alpha_i$ given by \eqref{e:nabla giq aalph-fne violation} for the steps $t_j$ with $j> s$ satisfying
\eqref{e:nabla gi aalph-fne step}.  Note that the steps $t_j$ for the nonsmooth function are only
implicity involved here through the constant $\tau_j$ in \eqref{e:sd f submon}.
Here too, the larger the step $t_j$, the larger the constant $\tau_j$, which leads to
a greater violation in \eqref{e:nabla giqr aalph-fne violation}.

\end{eg}

\subsection{Metric subregularity of the invariant Markov transport discrepancy}
\label{s:msr mtd}
Recall the inverse mapping
$\Psi^{-1}(y)\equiv \set{\mu}{\Psi(\mu)=y}$, which
clearly can be set-valued.  It is important to keep in mind
that an invariant measure need not correspond to a fixed point of any
individual mapping $T_i$, unless these have common fixed points.  See
\cite{HerLukStu23a, HerLukStu23b} instances of this.  Condition \eqref{e:Psi msr} together
with positivity of $\Psi$, that is $\Psi$  takes the value $0$ at
$\mu$ if and only if $\mu\in\inv\Pcal$, is the assumption that $\Psi$ is
{\em gauge metrically subregular for $0$ relative
to $\mathscr{P}_2(G)$ on $\mathscr{P}_2(G)$} \cite{Ioffe13}:
\begin{equation}\label{e:metricregularity}
d_{W_{2}}(\mu,\inv\Pcal) = d_{W_{2}}(\mu,\Psi^{-1}(0))\leq \Gcal(\Psi(\mu))
\quad \forall \mu\in \mathscr{P}_2(G).
\end{equation}
The gauge of metric subregularity $\Gcal$  is constructed
implicitly from  $\map{\theta_{\tau,\epsilon}}{[0,\infty)}{[0,\infty)}$
with parameters $\tau>0$ and $\epsilon\geq 0$ satisfying \eqref{eq:theta_tau_eps} and \eqref{eq:gauge}
for $\tau>0$ fixed.  In Proposition \ref{t:msr convergence} the parameter $\epsilon$ is
exactly the expected violation in the update mappings $\Phi(x,i) = T_i(x)$ which, as in Example \ref{eg:sd aafne} and \ref{eg:fb aafne},
are a$\alpha$-fne
(see \eqref{eq:paafne i.e.'} and
\eqref{eq:paafne i.e.});
the parameter $\tau$ is directly computed from the expected constant $\alpha$ via $\tau= (1-\alpha)/\alpha$.

When the mappings $T_i$ have at least one common fixed point, then the
feasible set is nonempty:
\begin{equation}\label{e:C}
 C\equiv \set{x\in G}{\Pbb\paren{x\in T_\xi x}=1} = \bigcap_{i\in\{1,2,\dots,M_m\}}\Fix T_i\cap G\neq \emptyset.
\end{equation}
Moreover, the set of measures supported on the feasible set is also nonempty, that is
\begin{equation}\label{e:Ccal}
 \mathscr{C}\equiv \set{\pi\in\inv\Pcal}{\Supp(\pi)\subset C} \neq \emptyset
\end{equation}
since the delta distribution centered at a point $x^*\in C$ is in $\mathscr{C}$.
The stochastic fixed point problem in this case is a consistent stochastic feasibility problem.
Here the Markov transport discrepancy $\Psi$ defined by \eqref{eq:Psi} simplifies considerably.
Indeed, take any $\pi\in\mathscr{C}$.  Then for $y\in\Supp\pi$ we have
$y= T_\xi y$ (almost surely) and \eqref{eq:Psi} becomes
  \begin{eqnarray}
\Psi(\mu)&\equiv& \inf_{\pi\in\inv\Pcal}\inf_{\gamma\in \Ccal_*(\mu,\pi)}
\left(\int_{G\times G}
\mathbb{E}\left[\norm{(x - T_\xi x) - (y - T_\xi y)}^2\right]\ \gamma(dx, dy)\right)^{1/2}
\nonumber\\
&=& \inf_{\gamma\in \Ccal_*(\mu,\delta_y)}
\left(\int_{G\times G}
\mathbb{E}\left[\norm{x - T_\xi x}^2\right]\ \gamma(dx, dy)\right)^{1/2}
\nonumber\\
&=&
\left(\int_{G}
\mathbb{E}\left[\norm{x - T_\xi x}^2\right]\ \mu(dx)\right)^{1/2}.
\label{eq:Psi_StoPBForBS}
 \end{eqnarray}
Clearly $\mu\in\mathscr{C}$ implies that $\Psi(\mu)=0$.
Condition \eqref{e:Psi msr} then reduces to
  \begin{eqnarray}
    d_{W_2}(\mu, \inv\Pcal)    &\le&
     \Gcal\paren{\left(\int_{G}
\mathbb{E}\left[\norm{x - T_\xi x}^2\right]\ \mu(dx)\right)^{1/2}}.
\label{e:Psi msr2}
  \end{eqnarray}
Suppose now that
$\mu$ is such that $d_{W_2}(\mu, \inv\Pcal)$ is attained at some $\pi\in \mathscr{C}$.
Then
  \begin{eqnarray}
    d_{W_2}(\mu, \mathscr{C}) = d_{W_2}(\mu, \inv\Pcal)    &\le&
     \Gcal\paren{\left(\int_{G}
\mathbb{E}\left[\norm{x - T_\xi x}^2\right]\ \mu(dx)\right)^{1/2}}.\label{e:Psi msr3}
  \end{eqnarray}
Writing this pointwise
(i.e., for $\mu=\delta_x$) reduces \eqref{e:Psi msr} to a recognizable (nonlinear) error bound,
albeit in expectation:
\begin{eqnarray}
d(x, C)  = \inf_{z\in C}\|x-z\| &\leq&
\Gcal\paren{\mathbb{E}\left[\norm{x - T_\xi x}\right]}\quad\forall x\in G.
\label{e:almost metricregularity StoPBForBS}
\end{eqnarray}

\section{Case Study:  likelihood maximization and X-FEL imaging}\label{s:em and rfi}
We specialize the composite minimization model to minimizing the negative log-likelihood function:
\begin{equation}\label{e:em}
 \ucmin{E(x)\equiv \sum_{j=1}^M-\log\paren{f\paren{y_j~;~\phi(x)}}}{x\in C_0}.
\end{equation}
Here $\map{\phi}{\Rn}{\Rm}$ is a nonlinear but smooth map modelling a proposed
probability distribution parameterized by $x\in \Rn$,  $C_0\subset\Rn$ is an
abstract constraint set,
$y_j\in \Rm$ ($j=1,2,\dots,M$) is an outcome, and $f\paren{y_j~;~\phi(x)}$
is the {\em likelihood} of observing the outcome $y_j$ with the probability distribution $\phi(x)$,
\begin{equation}\label{e:likelihood}
 f\paren{y_j~;~\phi(x)}\equiv \Pbb_{\phi(x)}\paren{Y_j=y_j}.
\end{equation}
The constraint $x\in C_0$ is formally handled by adding to the sum \eqref{e:em} the indicator function
\begin{equation}\label{e:g0 C0}
 g_0(x)\equiv \iota_{C_0}(x);
\end{equation}
the other functions in the objective \eqref{e:comp opt} are apparently
\begin{equation}\label{e:log-like}
 g_j(x)\equiv -\log(f(y_j; \phi(x))).
\end{equation}

The random variable $Y_j$ has some unknown true probability distribution.  The outcome $y_j$
is the result of a sample from the probability distribution.
We do not presume to approximate
the true probability distribution, although, obviously this is the whole point of problem
\eqref{e:em}.  Instead, we focus on numerical methods for solving \eqref{e:em} when the collection
of outcomes $\{y_1, y_2, \dots, y_M\}$ is too large to be handled all at once.
Interpretation of solutions to \eqref{e:em} as well as the limit of these solutions as $M\to\infty$ is
a topic for mathematical statistics.

As long as $f\paren{y_j~;~\phi(x)}$ is sufficiently smooth as a function of $x$ and the constraint set $C_0$ is
sufficiently regular (for instance, closed and convex), then the regularity assumption \eqref{eq:paafne i.e.} of
Proposition \ref{t:msr convergence} holds by the discussion in Examples \ref{eg:sd aafne} and \ref{eg:fb aafne};
hence the random function iteration Algorithm \ref{algo:rfi} is a reasonable candidate for mappings $T_i$ given by
\eqref{e:fbi}.
We will show this is indeed the case for the concrete application of X-FEL imaging.

\subsection{Stochastic tomography and single-shot X-FEL imaging}\label{s:xfel as rfi}
The goal of single-shot X-FEL imaging is to determine the electron density of a molecule from experimental
samples of its scattering probability distribution (i.e. diffraction pattern) \cite{Chap11,Bou12,ArdGru}. A
two-dimensional measurement device, partitioned into pixels, counts the occurrence of (coherently scattered)
photons in a three-dimensional domain in the far ﬁeld of a molecule that has been illuminated by a
short X-FEL pulse. The molecule
under observation is at a random and unobservable orientation relative to the measurement
device. The illuminating pulse is only long enough to cause a few (between 3 and 100) scattering
events from the interaction of the X-ray beam with the electrons of the molecule. The experiment is repeated
about $10^9$ times, each time with the molecule at a different random orientation.

Image reconstruction from single-shot X-FEL experiments is a stochastic
tomography problem with a nonlinear model for the data.  The stochastic aspect comes
not from the randomness of the molecule orientations (although this makes the problem harder),
but rather from the random samples of the scattered field that constitute the data. The
tomographic aspect comes from the challenge of reconstructing a three-dimensional object
from two-dimensional data.  Computed tomography with random orientations is not new
(see \cite{Bresler2000b} and  \cite{Bresler2000a}).  When the orientations are not known, or in the
case of single-shot X-FEL imaging unobservable, some earlier successful approaches applied to
Radon inversion (a linear imaging model) involved estimating the orientation first, and then
inverting the imaging model \cite{Panaretos09} and  \cite{SingerWu13}.

The problem investigated in
\cite{GruSchu23} differs from \cite{Bresler2000a, Bresler2000b, Panaretos09, SingerWu13} in three respects.
First and foremost, the model for the data in X-FEL imaging is nonlinear and
nonconvex: Fraunhofer diffraction models the field propagation
(mathematically equivalent to Fourier transformation, and thus linear), however
the measurements are only {\em samples} of the {\em intensities} of the (complex-valued) field
\cite{Born} and \cite{ArdGru20}. This leads to the second, equally important distinction, namely that
the model data is truly a low-count sample from an unknown distribution (i.e. field intensities),
which prohibits directly applying the already known techniques for
phase retrieval from coherent diffraction data, as is done in X-ray crystallography \cite{RobSal20}.
The third departure from previous approaches is that
the step of estimating the unknown orientation is avoided by computing the likelihood of the
observation at {\em all combinations of orientations}.

Considering only {\em coherently scattered photons}, the problem is to determine the true three-dimensional
electron density $\rho_*$ from a collection
of two-dimensional images $\{y_{s_1}, y_{s_2}, \dots, y_{s_M}\}$, where each $y_{s_j}$ is an outcome,
or sample from a two-dimensional cross section of the true three-dimensional probability distribution
in the far field $\phi_*$ at orientation $s_j\in SO(3)$, denoted $\phi_*^{s_j}$.
The continuous model for the probability distribution in the far field $\map{\phi_*}{\Rth}{\Rp}$ is derived from the
planar Fraunhofer scattering model \cite{Born}
\begin{equation}\label{e:true Fraun}
  \phi_*^s(k)\equiv |\ecklam{\Fcal\paren{\rho_*\circ R(s)}}(k)|^2, \quad k\in\Ecal.
\end{equation}
Here $\map{\rho_*}{\Rth}{\Rp}$ is the electron density, $\map{R(s)}{\Rth}{\Rth}$ is a rotation operator that
rotates the domain $\Rth$ by $s\in SO(3)$, $\map{\Fcal}{\Cbb^3}{\Cbb^2}$ is the Fraunhofer transform
modelling the propagation of the electromagnetic field to a two-dimensional surface in the far field, and,
in a nod to the physicists, we denote by
$k\in \Rth|_{\Ecal}$ a point on a two-dimensional surface $\Ecal$ (the Ewald sphere) in the far field.  Up to
scaling, $\Fcal$ is just the Fourier transform.  The function $\phi_*^s(k)$ is the {\em intensity} of
the field over long time spans and yields the probability of observing a photon at position $k$ when
the electron density $\rho_*$ is rotated by $s\in SO(3)$ and illuminated by an X-ray.

We approximate the true electron density with a function $\rho(\cdot; x)$ parameterized by
$x\in C_0\subset\Rn$.  The parameters are chosen so as to maximize the likelihood that random variables
$Y_{s_j}$ drawn from the true distribution $\phi_*$ match the observed outcomes $y_{s_j}$, or equivalently,
the optimal parameters $x^*$ minimize the negative log-likelihood functional \eqref{e:em} over some constraint set $C_0$.
This yields the following discretized model for the intensity
\begin{equation}\label{e:Fraun}
  \phi^s(k ; x)\equiv |\ecklam{\Fcal\paren{\rho(\cdot; x)\circ R(s)}}(k)|^2, \quad k\in\Rtw, x\in\Rn.
\end{equation}

In coherent diffraction imaging, the observations are over long time frames and the illuminating field
is assumed not to damage the object being observed (whether electron densities or cells is just a matter
of scale) \cite{RobSal20}.  In this scenario
one simply takes the observations $\{y_{s_1},y_{s_2},\dots,y_{s_M}\}$ to be identical to the probability
$\phi^s(k ; x)$ (when the rotation $s$ is known).  If enough photons are counted in each image, then the problem
reduces to the conventional optical phase retrieval problem for which successful numerical algorithms abound
\cite{LukSabTeb19}.  In a sense, this is the situation when each image $y_{s_j}$ is the true scattering
density at orientation $s_j$, in which case expectation maximization
really does find a best approximation of a reasonable parameterization to the true probability distribution.

In X-FEL imaging, the illuminating pulse destroys the molecule under observation, and so the time
is kept very short so that one can record scattering events (anywhere from $3$ to $100$) from the
molecule before it disintegrates.   The observations $\{y_1,y_2,\dots,y_M\}$ are therefore
simply counts of photons in a plane,  $y_j(k)\in \Nbb$ for $k\in \Rtw$. This is modeled as
a {\em Poisson point process}:
\begin{subequations}
\begin{eqnarray}\label{e:point process}
f(y_j(k); \phi^s(k; x)) =  \Pbb_{\phi^s(k; x)}\paren{Y^s_j(k)=y_j(k)}&\equiv& \frac{\phi^s(k ; x)^{y_j(k)}\exp\paren{-\phi^s(k ; x)}}{y_j(k)!}
\end{eqnarray}
where $f(y_j(k); \phi^s(k; x))$ is the likelihood function in \eqref{e:em}.
We have removed the dependence of the observation $y_j$ on the rotation $s$ since this is unobservable.
We will come to this in a moment.  Before accounting for the rotation, note that a single observation $y_j$
consists of all counts at locations $k\in \Rtw$, so the likelihood of $y_j$ with simultaneous counts at all
locations, assuming independence, is
\begin{eqnarray}
f(y_j; \phi^s(\cdot; x)) &=&  \prod_{k\in\Supp(y_j)}\frac{\phi^s(k ; x)^{y_j(k)}\exp\paren{-\phi^s(k ; x)}}{y_j(k)!}\nonumber\\
&=&\exp\paren{-\sum_{k\in\Supp(y_j)}\phi^s(k ; x)} \prod_{k\in\Supp(y_j)}\frac{\phi^s(k ; x)^{y_j(k)}}{y_j(k)!}.
\label{e:all point process}
\end{eqnarray}

\begin{remark}[modelling detector arrays]\label{r:pixels} Our numerical simulation adds the
discretization of the two-dimensional detector in the Fourier domain (k-space, for
the physicists) into pixels to the model \eqref{e:all point process}.
Within a single pixel the exact location of the photon cannot be determined.  Mathematically, this is
modelled by partitioning each observation $\Supp(y_j)$ into a rectangular grid.  To avoid clutter, we will not explicitly
represent this discretization in our derivation, but we will return to this source of error in our discussion of
the numerical results in Section \ref{s:results}.
\end{remark}

Returning to the issue of the rotations, there are really only two possibilities:  either one tries to
maximize the likelihood function with respect to $x$ {\em and} $s$, or one maximizes the likelihood with
respect to $x$, averaged over all possible rotations. The approach taken in \cite{GruSchu23} is the latter, namely
one considers
\begin{eqnarray}\label{e:av all point process}
f(y_j; \phi(\cdot; x)) \equiv \int_{SO(3)}f(y_j; \phi^s(\cdot; x)) \,d\mu(s).
\end{eqnarray}
\end{subequations}
Here we use the mean weighted by $\mu(s)$ to allow for the possibility of weighting
the rotations.  In the absence of any knowledge about
the scattering probability $\phi^s(\cdot; x)$ one would choose a uniform distribution
for the rotations.  But if one had a solution $x^*$ to \eqref{e:em}, this could be folded
back into \eqref{e:av all point process} to compute the average over the likelihood of
the rotations conditioned on $x^*$, effectively assigning a probability distribution to
the rotations for each outcome $y_j$.  Our demonstrations use only a uniform distribution
of directions on SO(3).

\begin{subequations}
An elementary calculation yields
\begin{eqnarray}
\nabla f(y_j; \phi(\cdot; x)) & = & \paren{\int_{SO(3)}\frac{\partial}{\partial x_i} f(y_j; \phi^s(\cdot; x)) \,ds}_{i=1,\dots,n}\nonumber\\
&=&
\paren{\int_{SO(3)}f(y_j; \phi^s(\cdot; x))\sum_{k\in\Supp(y_j)}
\paren{\frac{y_j(k)}{\phi^s(k; x)}-1}\frac{\partial}{\partial x_i}\phi^s(k; x) \,ds}_{i=1,\dots,n}.
\label{e: nabla av all point process}
\end{eqnarray}
Another elementary formal manipulation shows that
\begin{equation}\label{e:nabla Fraun}
\tfrac{\partial}{\partial x_i} \phi^s(k ; x)\equiv 2\Re\ecklam{\overline{\ecklam{\Fcal\paren{\rho(\cdot; x)\circ R(s)}}}(k)\cdot
\tfrac{\partial}{\partial x_i}\paren{\Fcal\paren{\rho(\cdot; x)\circ R(s)}}(k)}, \quad k\in\Rtw, x\in\Rn.
\end{equation}
One last application of the chain rule yields
\begin{equation}\label{e:nabla rho}
\tfrac{\partial}{\partial x_i}\paren{\Fcal\paren{\rho(\cdot; x)\circ R(s)}}(k) =
\paren{R(S)^*\Fcal^*\tfrac{\partial}{\partial x_i}\rho(\cdot; x)}(k)
\quad k\in\Rtw, x\in\Rn.
\end{equation}
The parameterization of the electron density $\rho$ is specified in the next section.
\end{subequations}

We can now  summarize the expression for
the gradients of the smooth functions\\
$g_j(x)\equiv -\log\paren{f\paren{y_j~;~\phi(x)}}$:
\begin{equation}\label{e:nabla gj em}
  \nabla g_j(x) = \frac{-1}{f\paren{y_j~;~\phi(x)}}\nabla f\paren{y_j~;~\phi(x)}
\end{equation}
where $\nabla f\paren{y_j~;~\phi(x)}$ is given by \eqref{e: nabla av all point process} for
$\tfrac{\partial}{\partial x_i} \phi^s(k ; x)$ given by \eqref{e:nabla Fraun} and
$\tfrac{\partial}{\partial x_i}\paren{\Fcal\paren{\rho(\cdot; x)\circ R(s)}}(k)$ given by
\eqref{e:nabla rho}.
The resolvent of the subgradient of the constraint-written-as-indicator function  $g_0$ is just the projector onto the set $C_0$:
\begin{equation}\label{e:Jg0 em}
  J_{g_0, 1}(x)\equiv P_{C_0}(x).
\end{equation}
This is all that is needed to implement Algorithm \ref{algo:rfi} with $T_i$ defined by \eqref{e:fbi}.

\begin{remark}[convergence of Algorithm \ref{algo:rfi} for X-FEL imaging]\label{r:XFEL convergence}
Assumptions (a) and (b) of Proposition \ref{t:msr convergence} are easily seen to be satisfied for this application.
Indeed, any twice continuously
differentiable function satisfies conditions \eqref{e:g reg} of Example \ref{eg:sd aafne}.  Moreover, if the
constraint $C_0$ is convex, then condition \eqref{e:sd f submon} is satisfied with $\tau_0=0$.
Hence the mappings
$T_i$ defined by \eqref{e:fbi} for the application of X-FEL imaging using model \eqref{e:em}
with a smooth parameterization of the electron density $\rho$
are a$\alpha$-fne with violation $\epsilon_i$ and constant $\alpha_i$ given by
\eqref{e:nabla giqr aalph-fne} for the steps $t_j$ with $j> 0$ satisfying
\eqref{e:nabla gi aalph-fne step}.
If conditions \eqref{t:msr convergence c} and \eqref{t:msr convergence b} of Proposition \ref{t:msr convergence}
hold, then we are assured that Algorithm \eqref{algo:rfi} is locally convergent in distribution;
if the gauge of metric subregularity
is linear, then by Corollary \ref{t:msr convergence - linear} we can expect local linear convergence in the
Wasserstein metric.
The numerical experiments in Section \ref{s:results} support this conclusion.
\end{remark}

Characterization of the invariant distributions is a topic of future research, but intuitively it should be clear that the more
observations $\{y_1, \dots, y_M\}$ and constraints that can be brought into the objective function \eqref{e:em}, the better
chance that the set of fixed points of Algorithm \ref{algo:rfi}
will contain only meaningful critical points of \eqref{e:em}.  The reasoning is as follows: if there is no data, then any feasible
parameters $x\in C_0$ will solve \eqref{e:em}.  As the number of observations grows, so too grow the number of functions $g_j$ in the log-likelihood
objective, and hence the set of critical points can only shrink.  If the set of critical points of \eqref{e:em} consists only of globally
optimal solutions, the invariant measures of the randomized algorithm \eqref{algo:rfi} {\em should} be supported on these solutions.
It is an open question to determine exactly when and in what way this holds.

\subsection{A Numerical Study}\label{s:results}
\begin{center}
	\begin{figure}[!htp]
     	(a) \includegraphics[width=6cm, height=6cm]{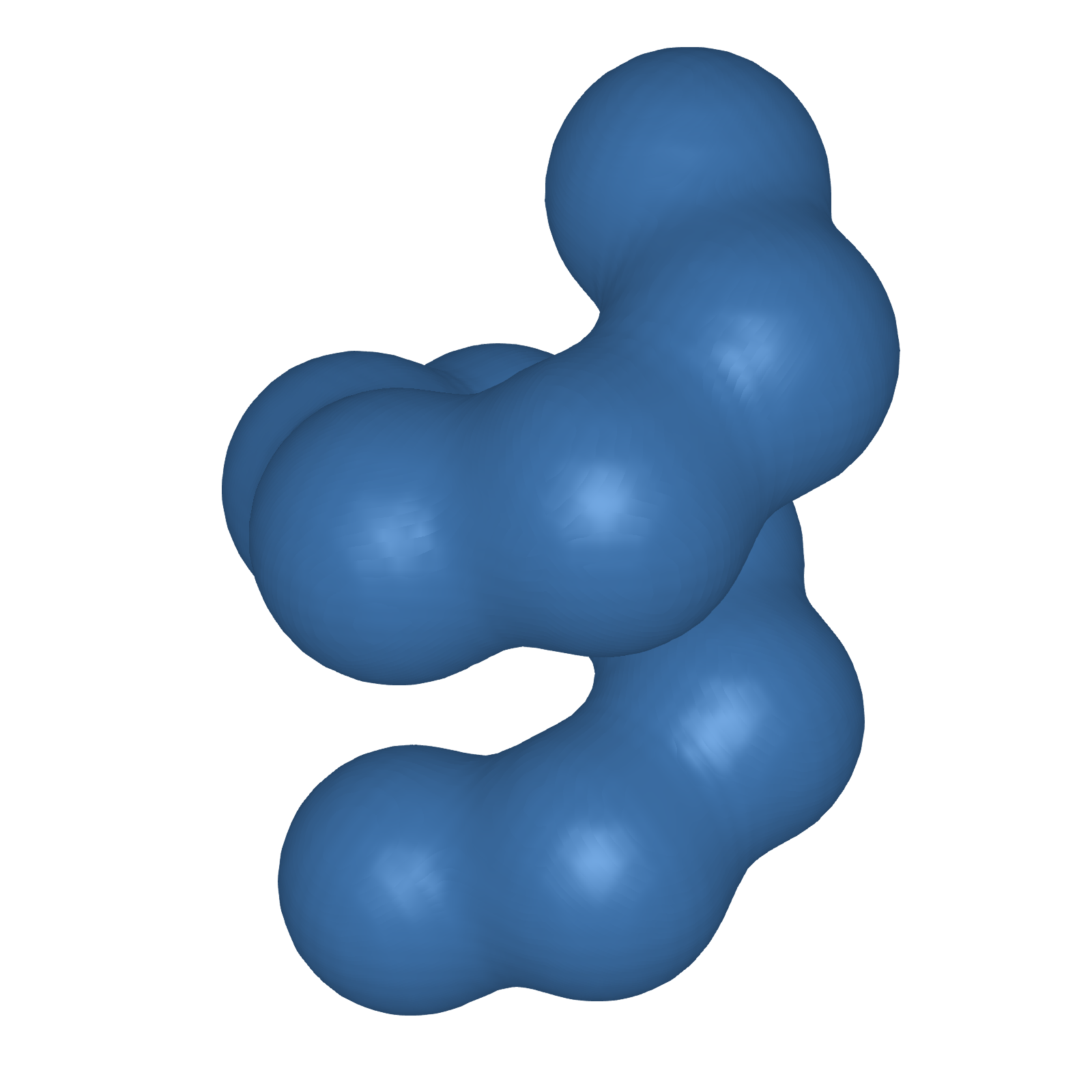}\hfill
     	(b) \includegraphics[width=6cm, height=6cm]{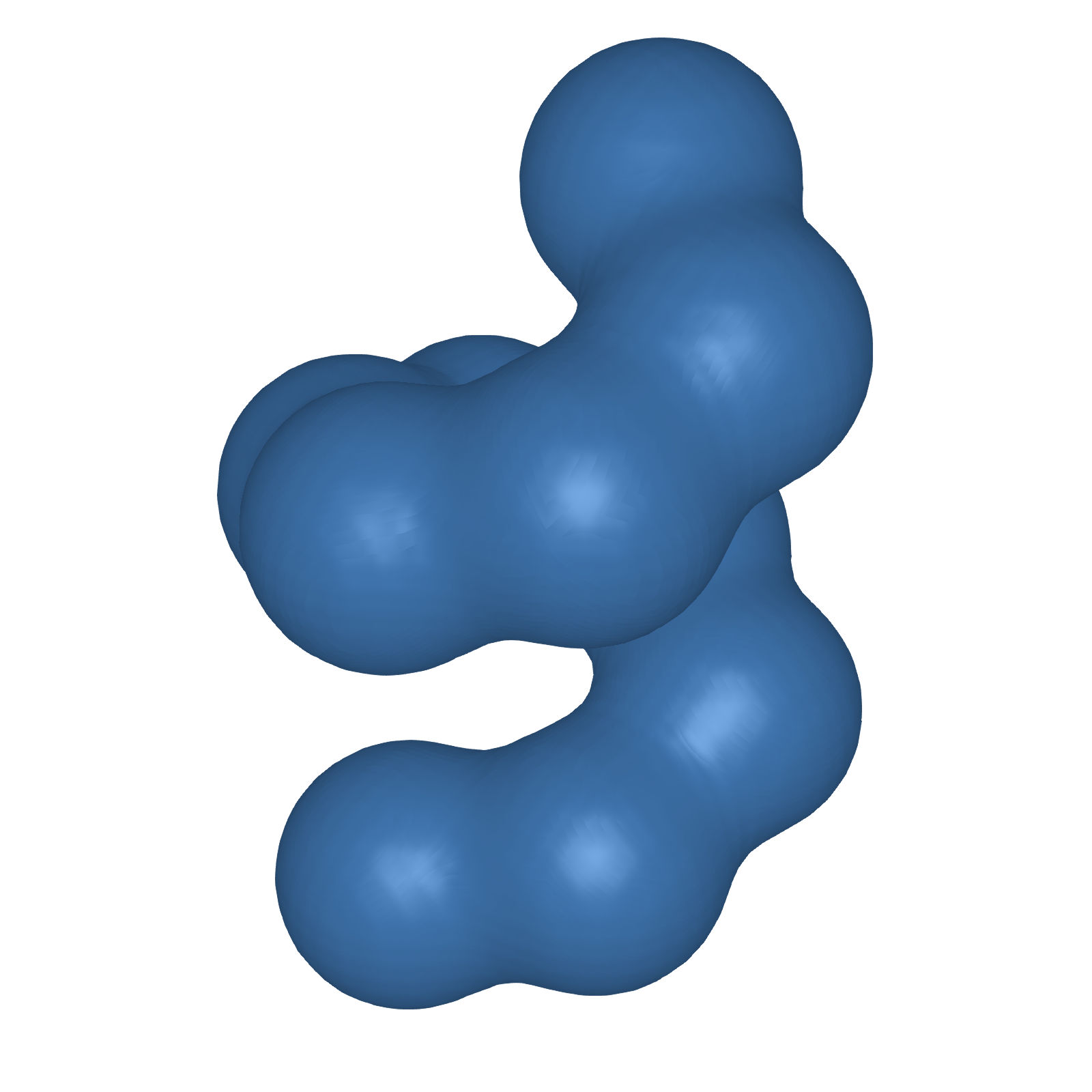}\\
     	(c) \includegraphics[width=6cm, height=6cm]{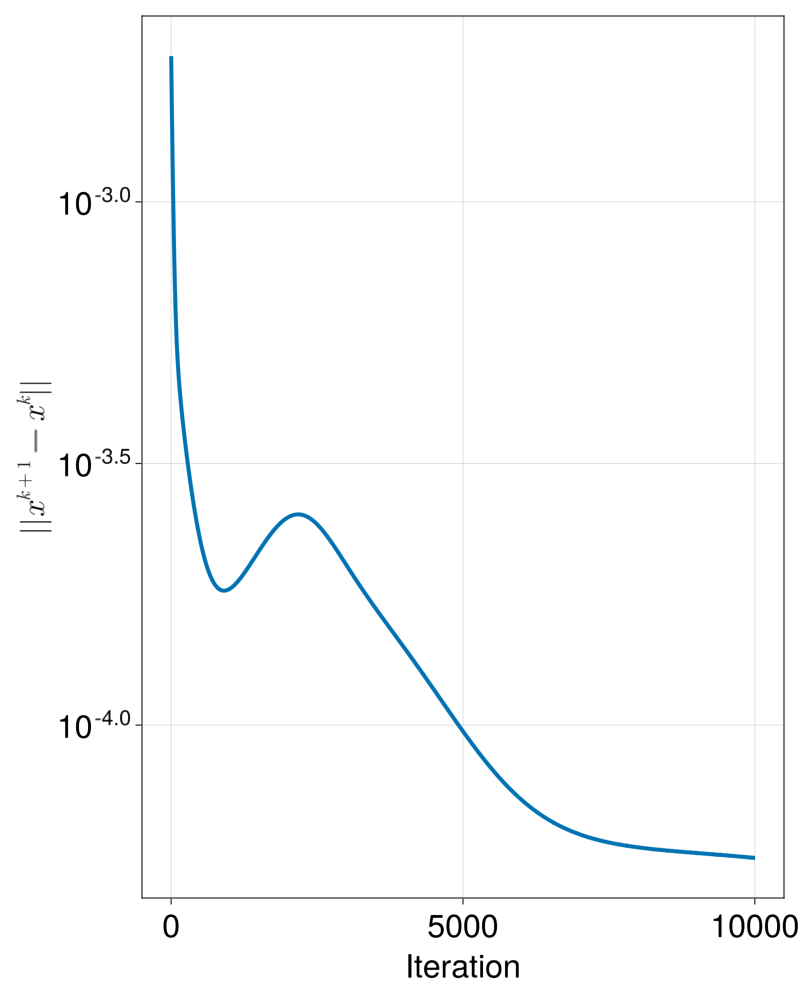}
		\caption{\label{fig:deterministic} Deterministic recovery.  Algorithm \ref{algo:rfi} with $m=M=10,000$, $t=0.1$.  (a) The true electron density.
		(b) The computed electron density at the last iterate using (a) as the starting point.
		(c) The iterate differences.}
	\end{figure}
\end{center}

A level surface of the electron density used to generate the data for our numerical experiment
is shown in \figr{fig:deterministic}(a).
The parameterized model for this is simply $10$ Gaussian balls where the parameters $x$ are the
centers of the balls:
\begin{equation}\label{e:rho}
 \rho(z; x)\equiv \sum_{i=1}^{10}\Bcal(z; x_i),\quad\mbox{ for }\quad
 \Bcal(z; x_i)\equiv \sigma*\exp\paren{-\tfrac{1}{2}\|z-x_{i}\|^2}
\end{equation}
%
where $x = \paren{x_1, x_2, \dots, x_{10}}$ with each $x_i = \paren{x_{i1}, x_{i2}, x_{i3}}\in\Rth$ and $\sigma = 0.519$.
Since the density has the form of Gaussians, the rotated field in the Fourier domain has the explicit formulation
\cite{Jain}
	\begin{subequations}\label{eq:electron density}
	\begin{eqnarray}\label{eq:rhohat}
	 \paren{\Fcal\paren{\rho(\cdot; x)\circ R(s)}}(k)&=& \sum_{i=1}^{10} \widehat{\Bcal}(k; \paren{R(s)x}_i)\\
	 \widehat{\Bcal}(k; \paren{R(s)x}_{i}) &=&
		\sqrt{2\pi}\sigma\exp\paren{-2\pi^2\|k\|^2}\exp\paren{-2\pi \im k\cdot (R(s)x_i)}
	 \label{eq:gaussians_hat}
	\end{eqnarray}
	\end{subequations}
	where $\im\equiv\sqrt{-1}$.
Just as a reminder, $k\in \Rth$ is a point in the Fourier domain restricted to a two-dimensional surface $\Ecal$ which
for simplicity we take to be a plane through the origin with a fixed orientation.
The gradient \eqref{e:nabla rho} then has the explicit representation
	\begin{eqnarray}
	 \tfrac{\partial}{\partial x_{ir}}\paren{\Fcal\paren{\rho(\cdot; x)\circ R(s)}}(k)&=&
       \tfrac{\partial}{\partial x_{ir}}\sum_{i'=1}^{10}\widehat{\Bcal}(k; \paren{R(s)x}_{i'})\nonumber\\
      &=&
	 -2\pi \im\sigma\paren{R(s)k}_r\widehat{\Bcal}(k; R(s)x_{ir})\quad (r=1,2,3).
	 \label{eq:nabla gaussians_hat}
	\end{eqnarray}
For this demonstration we do not introduce any constraints, so here $C_0=\Rn$, and the forward-backward
operator \eqref{e:fbi} reduces to just gradient descent \eqref{e:SDq}.

The numerical experiments were produced in the Julia programming language on 13th generation Intel\texttrademark
CORE\textregistered i9 CPUs together with Nividia GPUs.  Scripts for generating the numerical experiments
are available at \cite{PTb}.

We generated $M=10,000$ images/outcomes $y_j$ so that, on average, each image contained $15$
single scattering events across the measurement plane $\Ecal$.
\figr{fig:deterministic}(b)-(c) shows the result of applying Algorithm \ref{algo:rfi} to the
deterministic case: $m=M=10,000$
so that the index set $\Ibb=\{1,2,\dots,M\}$, the random variable $\xi_k$ takes only one value, and the
Markov operator reduces to the deterministic fixed point mapping $T\equiv \tfrac{1}{M}\sum_{j=1}^M\Id - t\nabla g_j$
where $\nabla g_j$ is given by \eqref{e:nabla gj em},  \eqref{e: nabla av all point process}, \eqref{e:nabla Fraun} and
\eqref{eq:nabla gaussians_hat}.  The mapping $T$ has fixed points at the critical points  of the objective
in \eqref{e:em} and these exist (assumption (a) of Proposition \ref{t:msr convergence}) by level-boundedness of the objective.
The invariant measures are supported exactly on these critical points:  $\Supp(\pi)\subset\Fix T$
for any $\pi\in \inv\Pcal$.  Assumption (b) of Proposition \ref{t:msr convergence} was shown to hold for all step sizes $t_j$ small
enough by the discussion in Remark \ref{r:XFEL convergence}.
The
Markov transport discrepancy $\Psi$ defined by \eqref{eq:Psi} reduces to
\eqref{eq:Psi_StoPBForBS}.  Convergence in distribution is actually pointwise deterministic convergence,
so from any initialization at a single point,  $\mu_0=\delta_{x^0}$, condition \eqref{e:Psi msr}
in Proposition \ref{t:msr convergence} reduces to \eqref{e:almost metricregularity StoPBForBS} which,
in this case, is
the usual deterministic error bound:
\begin{eqnarray}
d(x, \Fix T) &\leq&
\Gcal\paren{\norm{x - T x}}\quad\forall x\in G.
\label{e:error bound}
\end{eqnarray}
Under the assumption that this error bound holds for some gauge $\Gcal$, we expect, by Proposition \ref{t:msr convergence}
local convergence to a critical point of \eqref{e:em} for all fixed step lengths $t$ small enough.  The step
$t=0.1$ for these numerical experiments was found by trial-and-error.

The recovered image  in \figr{fig:deterministic}(b) was obtained from the initialization of the true electron
density shown in \figr{fig:deterministic}(a).  While it is difficult to see from the recovered object, the
iterate steps shown in \figr{fig:deterministic}(c) indicate that the numerical solution differs from the true object used to construct the model data.  There are several possibilities for this discrepancy
which will require further study.  One obvious source for this difference is that the information about the
true object is obtained only through finitely many random samples, so problem \eqref{e:em} is itself
a random instance, hence the recovered object, while obtained through a deterministic algorithm, is itself
a random variable.  The theory presented here can be applied to the limiting case of uncountably infinite
collections, i.e. in the limit as the sum in \eqref{e:em} or \eqref{e:comp opt} converges to an integral.
In fact, Assumption \ref{ass:1} and Proposition \ref{t:msr convergence} are based on
\cite[Assumption 2.1 and Theorem 2.6]{HerLukStu23b} which are formulated over possibly uncountably infinite
index sets.  Convergence of the laws of the iterates of Algorithm \ref{algo:rfi} for fixed sample sizes $m$ follows
by \cite[Theorem 2.6]{HerLukStu23b}.  The question as to whether the invariant distribution of the
Markov operator in this case corresponds to the true distribution $\rho_*$ has not been addressed here, but
we think that this should be true.  We will have more to say about this in the concluding remarks.
Another possible source of error is the operation of averaging
over all possible rotations in \eqref{e:av all point process}.  We do not compute the integral in
\eqref{e:av all point process}, but rather a discretized sum which introduces some unavoidable error.
An even more banal, but not uncommon source of error in practice is an incorrect value for one of the many
physical parameters used to construct the mapping $\Fcal$.  We are confident that with this synthetic
example, all such trivial model errors have been eliminated.

Some or all of these sources of error could be at play in practice with real data, but for the present
demonstration, these are beside the point.  We are focused only on finding solutions to model
\eqref{e:em} regardless of the
appropriateness of the chosen parameters, so long as bad parameter choices do not change the
mathematical properties of the objective (in particular smoothness).
In this context, \figr{fig:deterministic}(c) showing the differences between successive
iterates on a log-scale indicates local linear convergence to a critical point.  In the deterministic setting \cite[Theorem 2]{LukTebTha18} has
established that metric subregularity of $\Id-T$ for $0$ at fixed points
(that is, \eqref{e:error bound} in an appropriate form) is in fact {\em necessary} for local linear convergence.
\figr{fig:deterministic}(c) does not contradict the hypothesis that, in the deterministic case, the transport
discrepancy is metrically subregular with a linear rate, in which case the observed local linear convergence in
\figr{fig:deterministic}(c) is guaranteed by
Corollary \ref{t:msr convergence - linear}.  If the fixed point is, up to rotations, unique, then we can use the iterate
differences to obtain an estimate of the distance to the limit at termination.  Round-off  error for these
calculations starts to appear at about $10^{-5}$ since we are working on single-precision GPU's.
This deterministic calculation ($10,000$ iterations) took several minutes (less than 5).

\begin{center}
	\begin{figure}[!htp]
     	(a) \includegraphics[width=6cm, height=6cm]{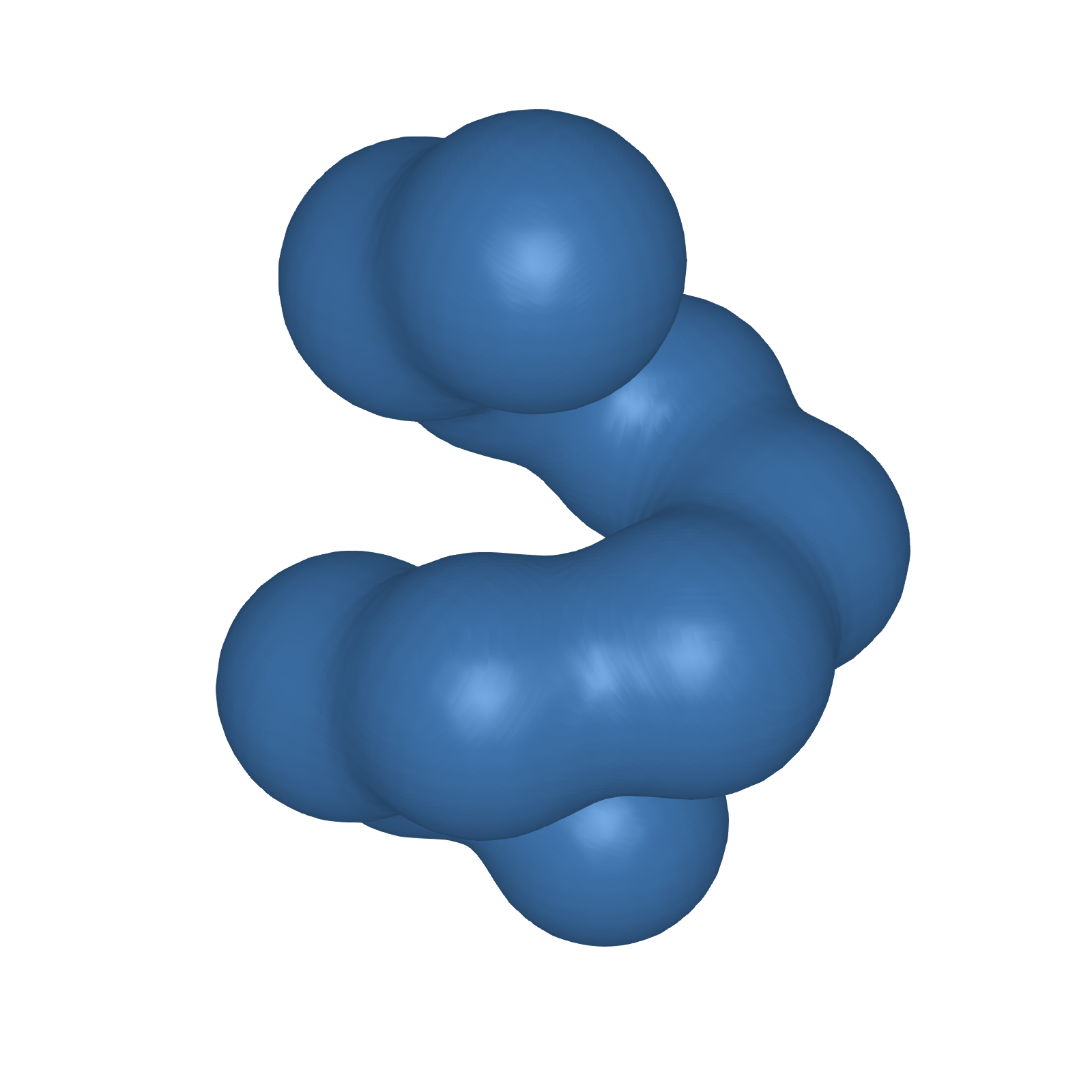}\hfill
     	(b) \includegraphics[angle=90, width=6cm, height=6cm]{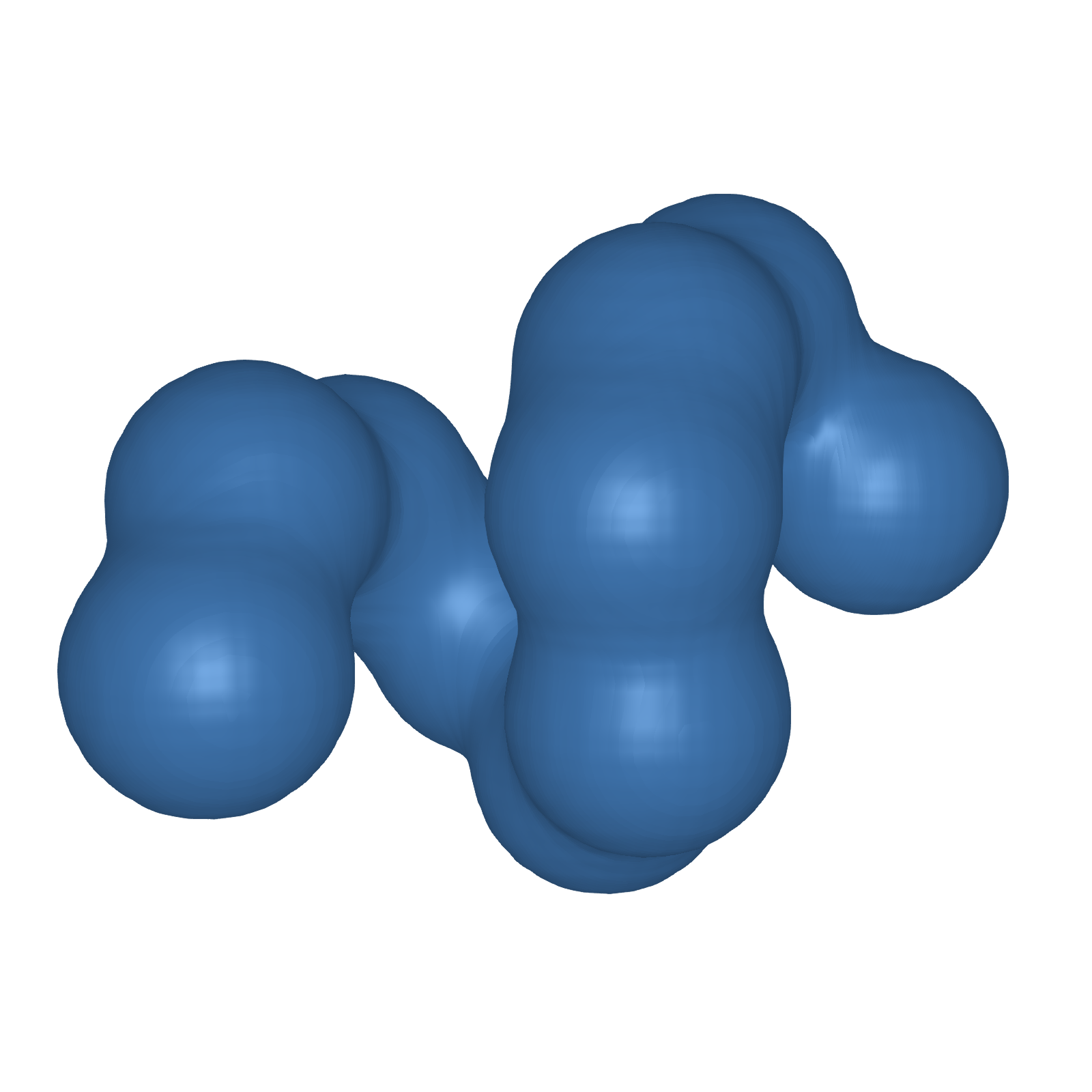}\\
     	(c) \includegraphics[width=6cm, height=6cm]{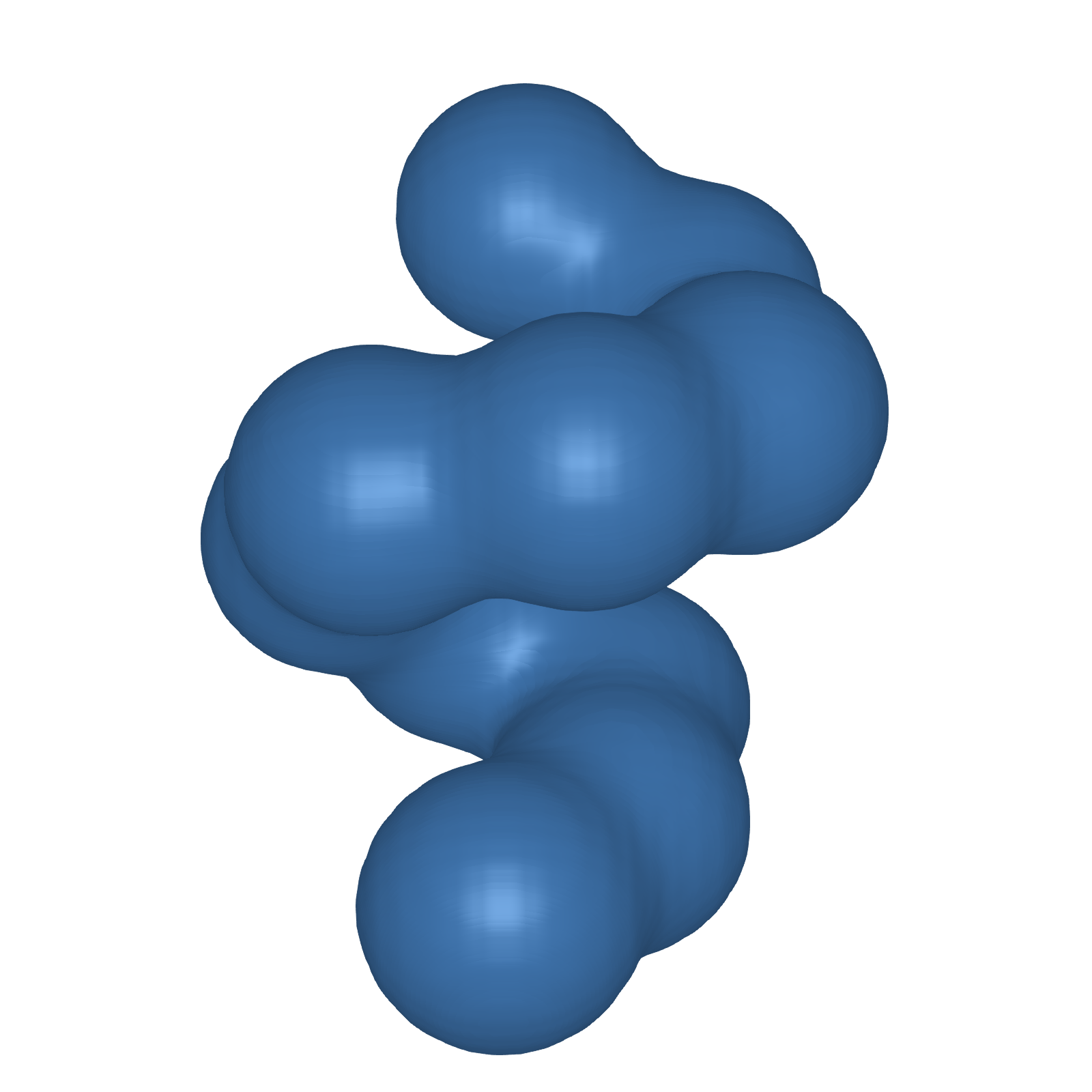}\hfill
        (d) \includegraphics[width=6cm, height=6cm]{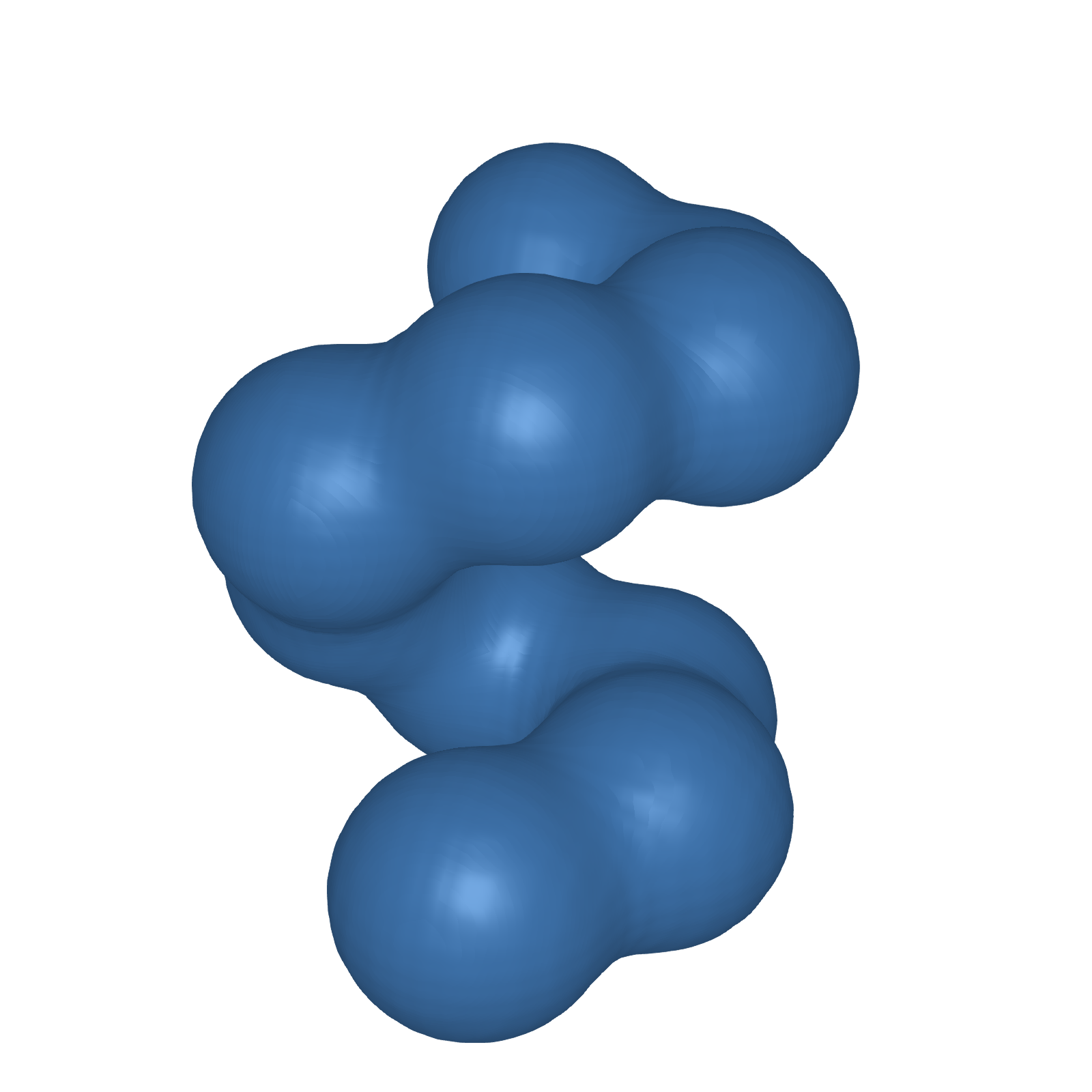}
		\caption{\label{fig:random recovery} Random recovery.  Algorithm \ref{algo:rfi} with $T_i$ given by \eqref{e:SDq} for $q=10$ and
		$t_j=0.1$ for all $j$.  The $M=10,000$ images are sampled with $|I_i|=m$ for  (a) $m=100$, (b) $m=500$, (c) $m=1000$, and (d) $m=5000$.
		Shown are the computed average electron densities at iteration $k=5000$.  The compute times
for the first $1000$ iterations ($1000$ outer iterations with $10$ inner iterations for each outer iteration) runs were:  $88$ seconds for $m=100$, $107$ seconds for $m=500$, $139$ seconds for $m=1000$, and
$418$ seconds for $m=5000$.  In comparison, the deterministic example with $10,000$ iterations and no sampling required
$494$ seconds}
	\end{figure}
\end{center}

For the randomized algorithm we sample, respectively,  $100$, $500$, $1000$, and $5000$ of the total $10,000$ images
($m=100, 500, 1000$ and $5000$) {\em without} replacement
and run $10$ iterations of gradient descent on each batch ($q=10$ in \eqref{e:SDq}).
The recovered density $\rho$ at the last (outer) iteration is shown in \figr{fig:random recovery}.
	\begin{figure}[!htp]
    \begin{center}
     	\includegraphics[width=6cm, height=6cm]{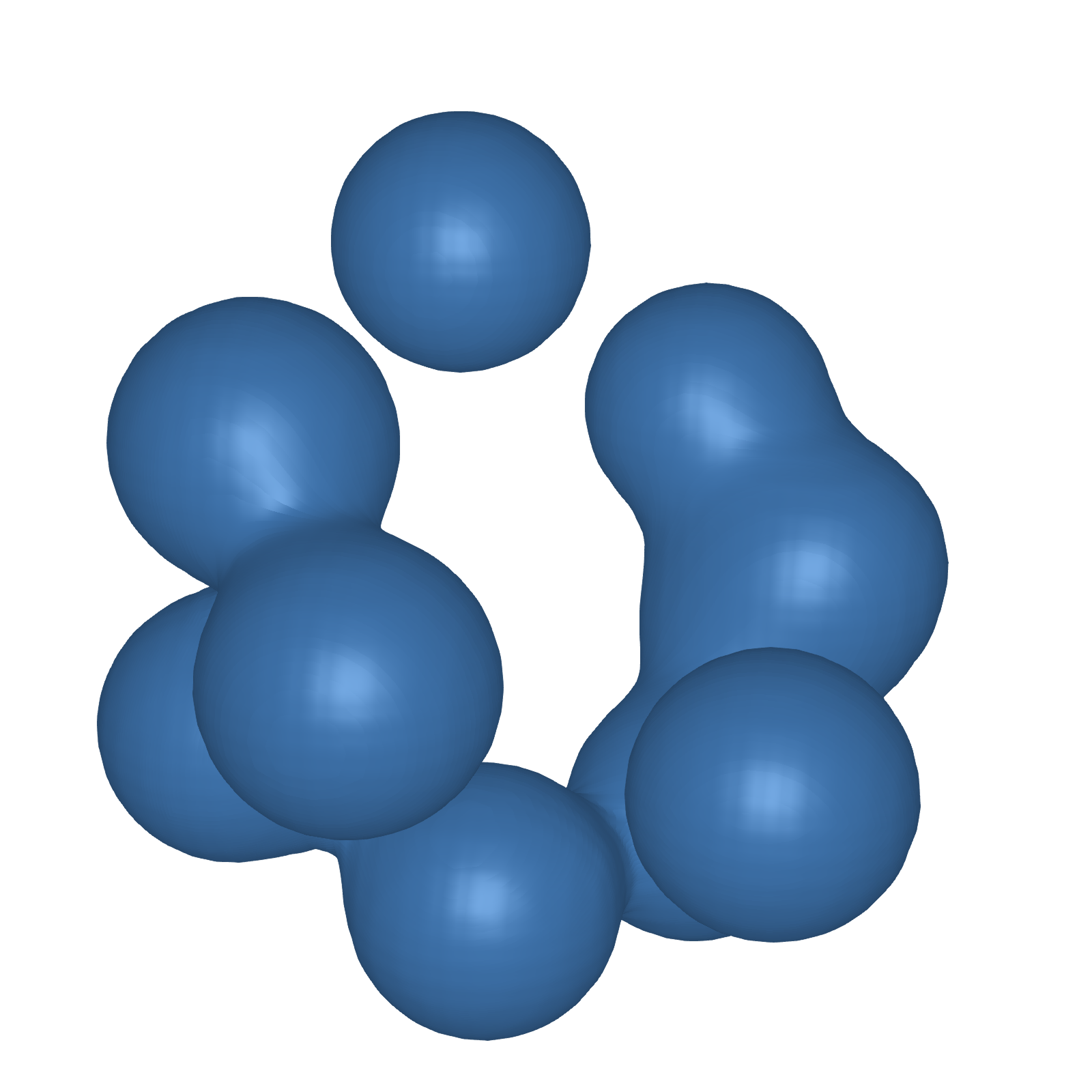}
		\caption{\label{fig:initialization} Random initialization.}
    \end{center}
	\end{figure}
Each numerical experiment was initialized by a randomly generated density, an example of which is shown in
\figr{fig:initialization}, and ran the same number of inner and outer iterations, namely $1000\times 10$;  this
was chosen to match the $10,000$ iterations for the deterministic example shown in \figr{fig:deterministic}.  The compute times
for these runs were:  $88$ seconds for $m=100$, $107$ seconds for $m=500$, $139$ seconds for $m=1000$, and
$418$ seconds for $m=5000$.  In comparison, the deterministic example with $10,000$ iterations and no sampling required
$494$ seconds.
\begin{center}
	\begin{figure}[!htp]
     	(a) \includegraphics[width=6cm, height=6cm]{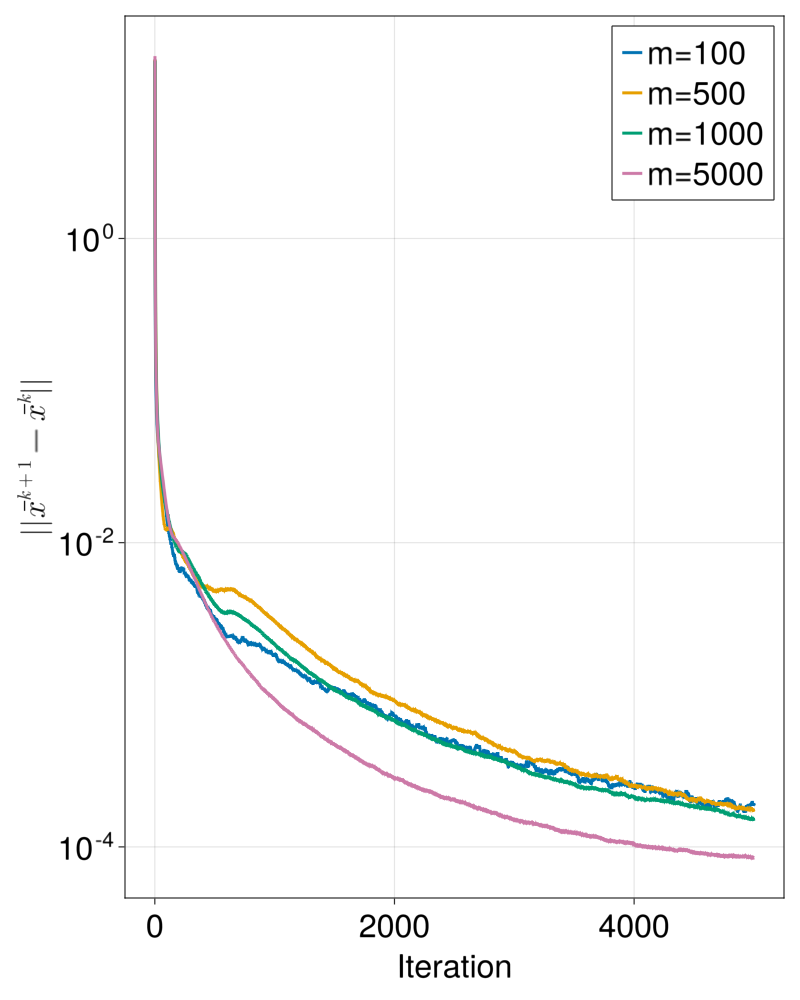}\hfill
     	(b) \includegraphics[width=6cm, height=6cm]{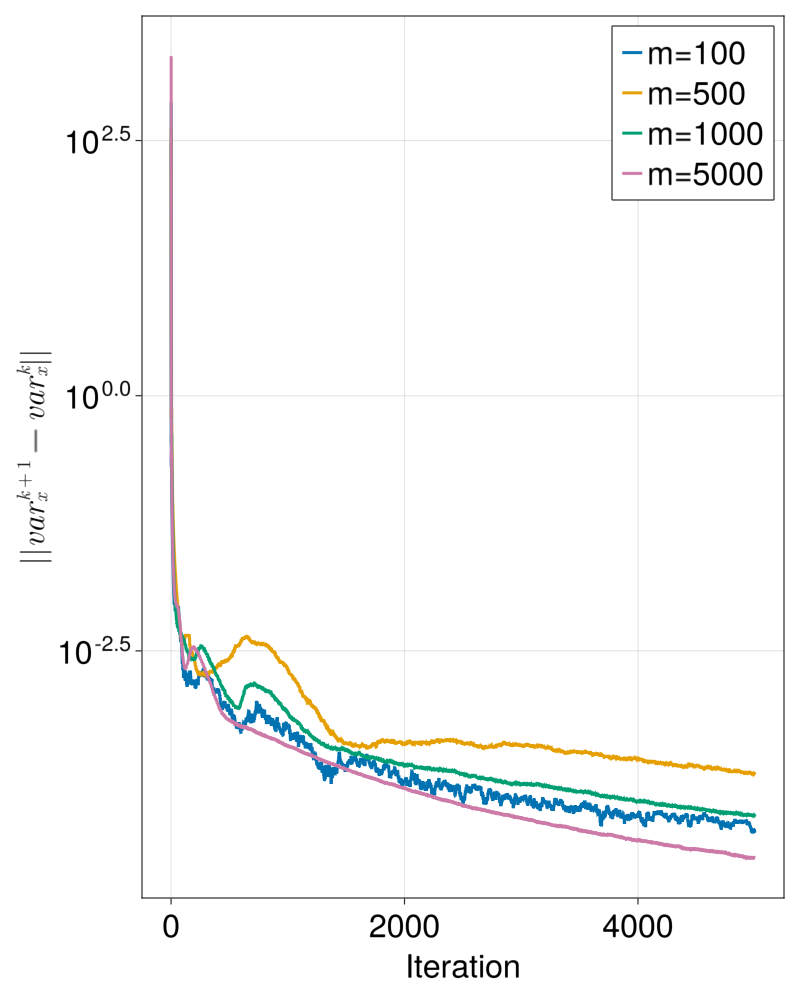}
		\caption{\label{fig:random convergence} Convergence behavior of the mean (a) and variance (b) of the iterates.
		Algorithm \ref{algo:rfi} with $T_i$ given by \eqref{e:SDq} for $q=10$ and
		$t_j=0.1$ for all $j$.  The $M=10,000$ images are sampled with $|I_i|=m$ for $m=100, 500, 1000$, and $5000$.}
	\end{figure}
\end{center}
As \figr{fig:random recovery} demonstrates, the reconstructions at iteration $5000$ are very similar, up to rotations and chirality,
to the deterministic reconstruction in \figr{fig:deterministic} at iteration $10,000$.
The Ces\`aro sum of the iterates and the variance of the iterate shown in \figr{fig:random convergence} was chosen to
indicate convergence of the first two moments of the laws of the iterates as predicted in Proposition \ref{t:msr convergence}
(or in the linear case \ref{t:msr convergence - linear}) since those results provide guarantees for convergence with respect
to the Wasserstein-2 probability measure, which yields rates of convergence for the first two moments of the distribution.

\section{Concluding remarks and open questions}
Several issues and directions for future research were mentioned above.  To conclude, we summarize these and add to this list.
In the context of X-FEL imaging, the directions for further study include: improving the numerical model without increasing
the computational complexity (too much);
 exploring more realistic parameterizations of the electron density;
 including constraints into the model \eqref{e:em}; and
 testing this approach on experimental data.  Regarding \ref{e:av all point process},
 it was mentioned that a solution $x^*$ can be used to assign a probability distribution to
the rotations of each outcome $y_j$, which can then be used for the weighted integral
in \eqref{e:av all point process}.
This is appealing both from the perspective of the application as well as posing
interesting mathematical challenges for the analysis of such an approach.

%

On the mathematical side, there are several difficult issues to tackle.  First and foremost is
either to develop efficient methods for computing the Wasserstein metric in the convergence
statement of Proposition \ref{t:msr convergence}, or to develop efficient approximations that
can be used to monitor algorithm performance.  Second but just as important is to characterize
the supports of the invariant measures of the Makov operator behind Algorithm \ref{algo:rfi} in
terms of the critical points of the determinisitc problem \ref{e:em} or, more generally \eqref{e:comp opt}.
This also has bearing on the analysis of the solutions to \eqref{e:em} in a stochastic context, in other
words, in the limit as the sum converges to an integral:  are the invariant measures of the Markov operator
unique up to rotation  and do these correspond to the true probability distribution $\rho_*$?
The dependence of the rates of convergence predicted in Proposition \ref{t:msr convergence} on the sample
size $m$ is also an important related issue.
Third, and more attainable in the near term is to verify the metric subregularity assumption of the
Markov transport discrepancy \eqref{e:Psi msr} by verifiable properties of its deterministic
counterpart.


\end{document}